\documentclass[11pt]{amsart}

\usepackage{amssymb}
\usepackage{amscd}
\usepackage{mathrsfs,eucal}

\usepackage{amsmath,hyperref,delarray,amsfonts}

\usepackage{times}

\def\nI{\relbar\hspace*{-0.4cm}\I\hspace*{0.1cm}}

\newcommand{\eop}{\hspace*{\fill}{\footnotesize $\blacksquare$}}

\def\<{\langle}
\def\>{\rangle}

\newcommand{\bt}{\begin{theorem}}
\newcommand{\et}{\end{theorem}}
\newcommand{\bc}{\begin{corollary}}
\newcommand{\bl}{\begin{lemma}}
\newcommand{\ec}{\end{corollary}}
\newcommand{\el}{\end{lemma}}
\newcommand{\bo}{\begin{observation}}
\newcommand{\eo}{\end{observation}}
\newcommand{\bp}{\begin{proposition}}
\newcommand{\ep}{\end{proposition}}
\newcommand{\br}{\begin{remark}}
\newcommand{\er}{\end{remark}}
\newcommand{\bmt}{\begin{maintheorem}}
\newcommand{\emt}{\end{maintheorem}}
\newcommand{\bq}{\begin{question}}
\newcommand{\eq}{\end{question}}

\newtheorem{theorem}{Theorem}[section]
\newtheorem{observation}[theorem]{Observation}
\newtheorem{corollary}[theorem]{Corollary}
\newtheorem{lemma}[theorem]{Lemma}
\newtheorem{proposition}[theorem]{Proposition}
\newtheorem{remark}[theorem]{Remark}
\newtheorem{maintheorem}{Main Theorem}
 
 \newtheorem{question}[theorem]{Question}

 \newcommand{\he}{\mathrm{id}}

\newcommand{\PG}{\ensuremath{\mathbf{PG}}}
\newcommand{\PGL}{\ensuremath{\mathbf{PGL}}}
\newcommand{\PSL}{\ensuremath{\mathbf{PSL}}}

\newcommand{\Aut}{\ensuremath{\mathrm{Aut}}}

\newcommand{\mS}{\ensuremath{\mathcal{S}}}
\newcommand{\mL}{\ensuremath{\mathcal{L}}}
\newcommand{\mO}{\ensuremath{\mathcal{O}}}

\newcommand{\mB}{\ensuremath{\mathcal{B}}}
\newcommand{\mK}{\ensuremath{\mathcal{K}}}
\newcommand{\mF}{\ensuremath{\mathcal{F}}}
\newcommand{\mQ}{\ensuremath{\mathbf{Q}}}

\newcommand{\mP}{\ensuremath{\mathcal{P}}}
\newcommand{\mG}{\ensuremath{\mathcal{G}}}

\newcommand{\mN}{\ensuremath{\mathcal{N}}}

\newcommand{\mM}{\mathcal{M}}
\newcommand{\mR}{\mathcal{R}}

\newcommand{\fix}{\mathrm{fix}}

\newcommand{\hW}{\ensuremath{\mathcal{W}}}
\newcommand{\hH}{\ensuremath{\mathcal{H}}}

\newcommand{\F}{\ensuremath{\mathbb{F}}}

\newcommand{\inc}{{\rm \tt I}}
\newcommand{\I}{\inc}

\newcommand{\id}{\ensuremath{\mathrm{id}}}

\newcommand{\hT}{\ensuremath{\mathbf{T}}}
\newcommand{\hF}{\ensuremath{\mathcal{F}}}

\newcommand{\proj}{\mbox{\rm proj}}

\title{Central aspects of skew translation quadrangles, I}


\author{Koen Thas}

\address{{Ghent University},
{Department of Mathematics},
{Krijgslaan 281, S25, B-9000 Ghent, Belgium}}
\email{koen.thas@gmail.com}
\thanks{}
\dedicatory{\texttt{Dedicated to my dear friend Ernie Shult}}
\date{}

\begin{document}

\maketitle

%\begin{abstract}
\subsection*{abstract}
Except for the Hermitian buildings $\hH(4,q^2)$, up to a combination of duality, translation duality or Payne integration, every known
finite building of type $\mathbb{B}_2$ satisfies a set of general synthetic properties, usually put together 
in the term ``skew translation generalized quadrangle'' (STGQ).
In this series of papers, we classify finite skew translation generalized quadrangles. In the first installment of the series, as corollaries
of the machinery we develop in the present paper,
\begin{itemize}
\item[|]
we obtain the  surprising result that any skew translation quadrangle of odd order $(s,s)$ is a symplectic quadrangle;\\
\item[|]
we determine all skew translation quadrangles with distinct elation groups (a problem posed by Payne in a less general setting);\\
\item[|]
we develop a structure theory for root-elations of skew translation quadrangles which will also be used in further parts, and which essentially 
tells us that a very general class of skew translation quadrangles admits the theoretical maximal number of root-elations for each member, and hence 
all members are ``central'' (the main property needed to control STGQs, as which will be shown throughout);\\
\item[|]
we solve the Main Parameter Conjecture  for a class of STGQs containing the class of the previous item, and which conjecturally coincides
with the class of all STGQs.\\
\end{itemize}
%\end{abstract}

\vspace*{1cm}
\tableofcontents

\newpage
\section{Introduction}

Let $L$ be a nilpotent $p$-group of class $2$, and suppose that $\Phi(L) = [L,L] \leq Z(L)$.
Define a map $\chi$ as follows:
\begin{equation}
\chi: V \times V \longrightarrow [L,L]: (a\Phi(L),b\Phi(L)) \longrightarrow [a,b], 
\end{equation}
where $V := L/\Phi(L)$ is seen as an $\mathbb{F}_p$-vector space.
Then $\chi$ is a bi-additive map. A fundamental question is to find minimal conditions  
on $L$ to conclude that $\chi$ is a bilinear {\em form} over $\mathbb{F}_q$, where $\vert \Phi(L) \vert = q$, a power of $p$. (So
$\chi$ then indicates how scalar multiplication works.) For instance, in \cite{Stroth}, Stroth formulated conditions on $2$-groups of class $2$ (which were
later corrected in \cite{Parker}) to force $\chi$ to be a form. Once one knows that $\chi$ is a bilinear form, one can use the geometry of 
$(V \times V,\chi)$ to study the structure of $L$. \\

For incidence geometry, such questions are very important. One of the most studied classes of buildings  is the class of so-called
``flock quadrangles'', and such quadrangles carry a subgroup $K$ of the automorphism group, for which (in the finite case)
$\vert K \vert = q^5$ for some prime power $q$, $\Phi(K) = [K,K] = Z(K)$ is elementary abelian of size $q$, and the map $\chi$ indeed is a (non-singular, alternating) bilinear form (``BAN-form''). We have that $K \cong \mathcal{H}_2(q)$ | the $5$-dimensional Heisenberg group over $\mathbb{F}_q$.
The group $K$ comes with two families $\mF$ and $\mF^*$ of subgroups which carry the geometric information of the quadrangle (i.e.,
the quadrangle can be  described in terms of cosets of elements of $\mF \cup \mF^*$). The elements of $\mF \cup \mF^*$ satisfy certain properties,
and vice versa, groups $K$ with families $\mF$ and $\mF^*$ of subgroups satisfying these properties act as certain automorphism groups of 
generalized quadrangles defined through this coset geometry construction. The pair $(\mF,\mF^*)$ is called a {\em Kantor family} in $K$.
In \cite{Isomflock}, I proved the following converse. 

\bt[\cite{Isomflock}]
\label{isomflock}
Suppose $H$ is a $p$-group of order $q^{5}$ for which $Z(H)  =  \Phi(H) =  [H,H]$ is elementary abelian of order $q$.
Suppose $H$ admits  a Kantor family of type $(q^2,q)$, and suppose $\chi$ defines a BAN-form over $\mathbb{F}_q$. Then $H \cong \mathcal{H}_2(q)$, and 
the corresponding generalized quadrangle $\mS$ of order $(q^2,q)$ is a flock quadrangle.
 \et
 
The main idea behind the proof is that the BAN-form produces a symplectic polar space (which is a higher rank Tits building) in the projective space coming from $V$, and  these geometries are well understood. Then going back to the quadrangle, one proves a characterizing property for flock quadrangles. \\
 
For a flock quadrangle (of order $(q^2,q)$), any defining Kantor family $(\mF,\mF^*)$ in $\hH_2(q)$ has the property that for $A^*, B^* \ne A^* \in \mF^*$, $A^* \cap B^* =: \mathbb{S}$ is a normal subgroup of $\hH_2(q)$ which is independent of the choice of $A^*, B^*$. Generalized quadrangles arising from Kantor families with this additional property are, by definition, {\em skew elation generalized quadrangles} (``STGQs''). This abstract class of quadrangles is very general:
except for the quadrangles associated to Hermitian varieties $\hH(4,q^2)$, every known finite generalized quadrangle is an STGQ 
up to a combination of duality, translation duality or Payne integration. 
 
So understanding the structure of STGQs (possibly up to the point of classification) is one of the main foundational challenges in the theory of rank two buildings.
And this is precisely the aim of the present series of papers.

\medskip
\subsection{From BN-Pairs to local Moufang conditions}

A famous result of Tits \cite{Titslect} states that an axiomatic spherical building of rank at least $3$ can be constructed through a group coset geometry construction  from a so-called ``(group with a) BN-pair'' (for the special case of axiomatic projective spaces, this boils down to the classical Veblen-Young theorem
which says that if the dimension is at least $3$, the space can be coordinatized over a skew field). In the other direction, each group with a BN-pair
gives rise to a building.  We mention that a group $G$ is said to have
a {\em BN-pair}\index{BN-pair} $(B,N)$\index{$(B,N)$}, where $B, N$ are subgroups of $G$, if the
following properties are satisfied:

\begin{itemize}
\item[(BN1)]
$\langle B,N \rangle = G$;
\item[(BN2)]
$H = B \cap N \lhd N$ and $N/H = W$ is a Coxeter group with distinct generators
$s_1,s_2,\ldots,s_n$;
\item[(BN3)]
$Bs_iBwB \subseteq BwB \cup Bs_iwB$ whenever $w \in W$ and $i\in\{1,2,\ldots,n\}$;
\item[(BN4)]
$s_iBs_i \ne B$ for all $i\in\{1,2,\ldots,n\}$.
\end{itemize}

The subgroup $B$, respectively $W$, is a {\em Borel subgroup}\index{Borel subgroup},
respectively the {\em Weyl group}\index{Weyl group}, of $G$. 
The natural number $n$ is called the \emph{rank}\index{rank} of the BN-pair; when the rank is $2$ the Weyl group $N/(B \cap N)$ is a
dihedral group of size $2m$ for some $m$. We say that the BN-pair $(B,N)$ is {\em of type $\mathbb{B}_2$}\index{BN-pair!of type $\mathbb{B}_2$} if $W$ is a dihedral group of size $8$.

From the data given in (BN1)--(BN4), one constructs maximal parabolics $P_i$, one for each generator $s_i$, and the element set of the associated 
building $\mB(G;B,N)$ consists of left cosets of the maximal parabolics, while such elements are incident if they have nontrivial intersection. 
Then $G$ acts naturally as an automorphism group on $\mB(G;B,N)$ by left multiplication. Sphericality expresses the fact that $W$ is a finite group (which translates geometrically in the fact that apartments are finite). One can prove that higher rank buildings can be 
described completely by their rank $2$ residues \cite{Titslect}, which are buildings of rank $2$, i.e., {\em generalized polygons} \cite{POL}.

In \cite{Titslect}, Tits obtains a splitness result for each BN-pair associated to a spherical building of rank at least three; interpreted on the rank $2$ residues, 
one obtains that the residual  polygons are {\em Moufang polygons}. Much later, the classification of Moufang polygons by Tits and Weiss \cite{TitsWeiss} became a flagship in incidence geometry, and showed to be more difficult than the classification of higher rank spherical BN-pairs (having Tits's classification  \cite{Titslect}
as a consequence).

The Moufang condition and local variations has become the central group theoretical condition in incidence geometry, as (e.g.) the well-known and much-studied notions of {\em translation plane} \cite{TP}  (planes that admit all local Moufang conditions at a line), {\em elation quadrangle} \cite{LEGQ} (the natural generalization to quadrangles, which enabled Payne and Thas to give a new proof of the classification of finite Moufang quadrangles \cite{PT,PTsec}), {\em translation quadrangle} \cite{TGQ} (another ``abelian'' variation), {\em Moufang polygon} \cite{TitsWeiss},
etc. show. The reader observes that in the rank $2$ case, one needs to impose extra (for instance group 
theoretical) conditions in order to have a grasp on these objects, much opposite to the higher rank case where large automorphism groups are always
present, as we have seen. The concept of free construction is just one of the many examples which underlines this view. 

Although in the theory of finite generalized polygons, local Moufang conditions have not produced new examples if the gonality is at least $6$, the 
two other cases (gonality $3$ and $4$) have proved to be conceptually different: in fact, recent work has shown that especially in gonality $4$ (i.e., for 
{\em generalized quadrangles}) automorphism groups can behave in a deep and unpredictable way. Even more, in the last twenty years, 
classes of examples have been (re-)discovered which satisfy a combination of local Moufang conditions, though not satisfying the global Moufang 
condition. Perhaps the main role has been played (in one way or another) by Payne's MSTGQs \cite{9} (see further), which is a very rich class
containing the important flock quadrangles, that is ``locally very Moufang''. All these examples are examples 
of STGQs, and by recent work of the author and also by work presented here, one observes that the class of STGQs strictly contains the MSTGQs. \\

In fact, this entire paper and its sequels should be seen as a new vista on local Moufang conditions in finite and infinite generalized quadrangle theory. 
(More precisely, we call an STGQ $\mS^x$ {\em central} if each symmetry about $x$ is contained in the center
of any elation group w.r.t. $x$.  We refer to ``Centrality conjecture'' as the conjecture that any STGQ is central.  
The Centrality conjecture is one the main actors in this series of papers | it will be obtained for (*)-STGQs (cf. below) in the present paper.
Centrality gives rise to a number of local Moufang conditions (and vice versa), which enable one to control this large class of quadrangles.)

\medskip
\subsection{The known results}

The first positive result towards STGQ classification is the one, independently obtained by Chen and  Hachenberger, on 
the parameters of an STGQ (solving a well-known conjecture of Payne):

\bt[Chen \cite{XC}, see also \cite{LEGQ}; Hachenberger \cite{Hach}]
\label{par}
The parameters $s,t$ of an STGQ are powers of the same prime.
\et

A more recent one, by Bamberg, Penttila and Schneider, is, by my knowledge, a first
``complete'' classification result, after assuming conditions on the parameters:

\bt[Bamberg, Penttila and Schneider \cite{BAM}]
An STGQ of order $(s,p)$, with $p$ a prime, is either classical or a flock quadrangle.
\et

In fact, in \cite{BAM} the same conclusion is proved for EGQs, but it follows easily from the assumptions that the quadrangle is an STGQ.
The fact that $t = p$ is a prime forces that either $s = t$ or $s = t^2$. The first case is known to lead to a quadric; in the second case
the authors show that the aforementioned bi-additive map is a BAN-form. \\

The final classification result we know of (besides Theorem \ref{isomflock}), is the very elegant

\bt[Ghinelli \cite{Ghi}]
\label{Ghi}
A square STGQ with odd order and of symplectic type is classical.
\et

One of the results of the present paper tells us that square STGQs with odd order {\em always} are of symplectic type, 
thus finishing the classification of odd square order STGQs (we indicate an independent conclusion why the STGQ then is classical). 
I announced this result back in May 2007 (cf. below). In fact, {\em all} (*)-STGQs (except for even square ones) will be shown to be 
of symplectic type in a geometric sense.

\br[Quadrangular representations]
{\rm
We say that a group $K$ has a (faithful) {\em quadrangular representation} (of type $(s,t)$) if there is a generalized quadrangle $\mS$, a point $x$ of $\mS$, and
an injection $\iota: K \longrightarrow \Aut(\mS)$ such that $(\mS^x,\iota(K))$ is an EGQ (of order $(s,t)$). The representation is {\em skew} if $(\mS^x,\iota(K))$ is an STGQ.
One of the goals in this series is to classify and understand (finite and infinite) groups which have a skew quadrangular representation, as well as
the isomorphism classes of possible representations of such groups. Another fundamental question that arises is which groups (if any)
can have quadrangular representations of different type. The only examples I know of are certain elementary abelian groups.
(Similar questions arise in the theory of Singer groups for projective spaces. In that context, 
existence of groups which act as Singer groups on finite projective spaces of different dimensions yield strange number theoretical identities.) \\
}
\er

\subsection{Structure of part I}

The paper is organized as follows.\\

\quad 0| In a first part (0) of the paper we introduce some elementary combinatorial features which will be used frequently throughout. By 
no means we are complete | we only define the language which is most often used, and refer to other sources on the way, if needed.
The concept of skew translation generalized quadrangle is described, as are
the known (classes of) examples of skew translation generalized quadrangles in some 
detail (we also mention the known elation groups which arise). These descriptions partly come from very recent
results, and contain new results on the local structure of one of these elation groups which sheds new light on the possible actions of skew elation groups.\\

\quad I| In a second part (I), we develop a basic structure theory for fixed points configurations for elements of skew elation groups.
We introduce a crucial property (called ``(*)''), which is enjoyed by most of the known STGQs. However, the aforementioned
local result shows that (surprisingly) there are examples which do not have (*). This obstruction motivates us to introduce ``generic STGQs'' later in the paper.
As a corollary of the fixed point theory, 
we show that once an STGQ has (*), we can deduce precise (and very strong, as the reader will see) information on the number
of fixed lines of elation group elements. This is one of the main tools for the next part.\\

\quad II| In that part (II), we discuss the existence of root-elations in STGQs, mostly in the case that (*) holds. The aforementioned fixed line 
theorem allows us to show that (*)-STGQs (assumed that when the STGQ is square, we are not in even characteristic) satisfy all local
Moufang conditions with respect to the elation point. We then deduce that such STGQs are central. We indicate several entirely different ways 
to do this, each with its own advantages.\\

\quad III| Next, we obtain one of the central results of this paper (and of the entire series): we show that all STGQs $(\mS^x,K)$ of odd square
order are isomorphic to $(\hW(t),\hH_1(t))$, that is, all come from a symplectic polar space in $3$ dimensions. So the only groups with a faithful skew 
quadrangular representation of type $(t,t)$, $t$ odd, are $3$-dimensional Heisenberg groups.
This generalizes Theorem \ref{Ghi} to a large extent. \\

\quad IV| In a fifth part, we consider STGQs which have ideal subGQs (through the elation point), and show that this assumption is enough 
to deduce strong structural information, including local Moufang conditions and centrality. We apply this theory to consider STGQs which admit
different elation groups (after a question of Payne), and solve his question by 
showing that such STGQs arise from Hermitian varieties. Also, further in the paper and in the series, the subGQ theory will be highly useful.\\

\quad V| Finally, we study the category $\mG$ of ``generic STGQs'' | STGQs which do not enjoy (*), but which we enhance with two very general natural properties
that generalize (*). (In fact, one of the main aims is to show that the category of generic STGQs coincides with that of non-(*) STGQs.) We show 
that every object of $\mG$ has ideal subGQs, so that the previous part can be applied. As a corollary, we show that for a (*)-STGQ or generic STGQ, either the possible parameters have the form $(t,t)$ or $(t^2,t)$ (with $t$ a 
prime power), or the STGQ is ``abelian'' (in a very specific case). This is a first big step which aims, eventually, at showing that the same 
arithmetic conclusion holds for any STGQ (Main Parameter Conjecture). \\

In an appendix, we make some initial remarks on ``abelian STGQs'' and a local version of (*), mainly preparing later work in the series.\\

\subsection{Acknowledgements}

I presented some of the main results of this paper as a series of (then partly conjectural) results in a lecture, at the conference ``Buildings and Groups'' (Ghent, May 2007), which was organized by 
Peter Abramenko, Bernhard M\"{u}hlherr and Hendrik Van Maldeghem.
Some of the ideas presented at that conference were developed when I 
was hosted, together with S. De Winter and E. E. Shult, by the Mathematisches Forschungsinstitut Oberwolfach (MFO) in the 
Research in Pairs program (April 2007).\\

\newpage
\section{Some preliminaries}
\label{STGQ}

We start this section with introducing some combinatorial and group theoretical notions.

\subsection{Elementary combinatorial preliminaries}

We tersely review some basic notions taken from the theory of  generalized quadrangles, for the sake of convenience.

\subsubsection{Finite generalized quadrangles}

Let $\Gamma$ be a thick generalized quadrangle (GQ). It is a rank $2$ geometry $\Gamma = (\mP,\mB,\I)$ (where we call the elements of $\mP$ ``points'' and those of $\mB$ ``lines'')
such that the following axioms are satisfied:
\begin{itemize}
\item[(a)]
there are no ordinary digons and  triangles contained in $\Gamma$;
\item[(b)]
each two elements of $\mP \cup \mB$ are contained in an ordinary quadrangle;
\item[(c)]
there exists an ordinary pentagon.
\end{itemize}
It can be shown that there are exist constants $s$ and $t$ such that each point is incident with $t + 1$ lines and each line is incident with $s + 1$ points. We say that $(s,t)$ is the {\em order} of $\Gamma$. \\

Note that an ordinary quadrangle is just a (``thin'') GQ of order $(1,1)$ | we call such a subgeometry also ``apartment'' (of $\Gamma$). If $s,t > 1$, $\mS$ is {\em thick}\index{thick}.\\

If $s = t$, then $\mS$ is also said to be {\em of order $s$}\index{order}; we call GQs of order $s$ for some $s$ {\em square generalized quadrangles}.\\

Suppose $p \nI L$. Then by $\proj_Lp$\index{$\proj_Lp$}, we denote the unique point on $L$ collinear with $p$.
Dually, $\proj_pL$\index{$\proj_pL$} is the unique line incident with $p$ concurrent with $L$.

\subsubsection{Roots and i-roots}

Let $A$ be an apartment of a GQ $\Gamma$. A {\em root} $\gamma$ of $A$ is a set of $5$ different elements $e_0,\ldots,e_4$ in $A$ such that $e_i \I e_{i + 1}$ (where the 
indices are taken in $\{0,1,2,3,4\}$), and $e_0, e_4$ are the {\em extremal elements} of $\gamma$. There are two types of roots, depending on whether the extremal elements are lines or points; in the second case we speak of {\em dual roots} to make a distinction between the types. Also, a (dual) root $\gamma$ without its extremal elements | the {\em interior} of $\gamma$ | is denoted by $\dot{\gamma}$ and called (dual) {\em i-root}.

\subsubsection{Point-line duality}

There is a {\em point-line duality}\index{point-line duality} for GQs of order $(s,t)$ for which in any
definition or theorem the words ``point'' and ``line'' are interchanged and
also the parameters. (If $\mS = (\mP,\mB,\I)$ is a GQ of order $(s,t)$, $\mS^D = (\mB,\mP,\I)$\index{$\mS^D$} is a GQ of order $(t,s)$.)

\subsubsection{Collinearity, concurrency and regularity}

Let $p$ and $q$ be (not necessarily distinct) points of the GQ $\mathcal{S}$; we write
$p \sim q$ and call these points {\em collinear}\index{collinear}, provided that there
is some line $L$ such that $p \I L \I q$. Dually, for $L, M \in \mB$, we write $L \sim M$\index{$L \sim M$}
when $L$ and $M$ are {\em concurrent}\index{concurrent}. 
For $p \in \mP$, put\index{$p^{\perp}$}
\begin{equation} p^{\perp} = \{ q \in \mP \vert q \sim p \}\end{equation}

(and note that $p \in p^{\perp}$). For a pair of distinct points $\{ p,q \}$, we denote $p^{\perp} \cap q^{\perp}$ also
by $\{ p,q \}^{\perp}$\index{${p,q}^{\perp}$}. Then $\vert \{ p,q \}^{\perp} \vert = s + 1$ or $t + 1$,
according as $p \sim q$ or $p \not \sim q$, respectively.
For $p \ne q$, we define
\begin{equation} \{ p,q \}^{\perp \perp} = \{ r \in \mP \vert r \in s^{\perp}\index{$p^{\perp\perp}$} \ \ \mbox{for all}\ \  s \in
\{ p,q \}^{\perp} \}.\end{equation}

When $p \not \sim q$, we have that $\vert \{ p,q
\}^{\perp\perp} \vert = s + 1$ or $\vert \{ p,q \}^{\perp\perp} \vert \leq t + 1$ according as $p \sim q$ or $p \not \sim q$, respectively. If $p \sim q$, $p \ne q$, or if $p \not \sim q$ and $\vert \{ p,q \}^{\perp\perp} \vert
= t + 1$, we say that the pair $\{p,q\}$ is {\em regular}\index{regular!pair}.
The point $p$ is {\em regular}\index{regular!point}\index{regular!line} provided $\{p,q\}$ is regular for every $q \in \mP \setminus \{ p \}$.
Regularity for lines is defined dually.
One easily proves that either $s = 1$ or $t \leq s$
if $\mathcal{S}$ has a regular pair of noncollinear points; see 1.3.6 of \cite{PT}.

%{\bf The GQ $W(q)$}.\quad
%Consider the $3$-dimensional projective space $\PG(3,q)$ over the finite field with $q$ elements $\mathbb{F}_q$. The points of $\PG(3,q)$ together with the totally isotropic lines with respect to a symplectic polarity form a GQ of order $q$, denoted $W(q)$.
%All the points of $W(q)$ are regular \cite{PT}, and this characterizes $W(q)$ as the only GQ of order $s = q \ne 1$ with this property.\\

\subsubsection{Antiregularity}

The pair of points $\{x,y\}$, $x \not\sim y$, is {\em antiregular}\index{antiregular!pair} if
$\vert \{x,y\}^{\perp} \cap z^{\perp} \vert \leq 2$
for all $z \in \mP \setminus \{x,y\}$. The point $x$ is {\em antiregular}\index{antiregular!point}
if $\{x,y\}$ is antiregular for each
$y \in \mP \setminus x^{\perp}$.

\subsection{Automorphism groups}

\subsubsection{Automorphisms and whorls}

An {\em automorphism}\index{automorphism} of a GQ $\mS = (\mP,\mB,\I)$ is a permutation of $\mP \cup \mB$ which
preserves $\mP$, $\mB$ and $\I$. The set of automorphisms of a GQ $\mS$ is a group, called
the {\em automorphism group}\index{automorphism!group} of $\mS$, which is denoted by $\Aut(\mS)$\index{$\Aut(\mS)$}.

A {\em whorl}\index{whorl} about a point $x$ is just an automorphism fixing it linewise.
 A point $x$ is a {\em center of transitivity}\index{center!of transitivity} provided that the group of whorls about $x$ is transitive
on the points of $\mP \setminus p^{\perp}$.

A {\em symmetry} with {\em center} $x$ is a whorl about $x$ which fixes $x^{\perp}$ pointwise. If the GQ $\mS$ is finite of 
order $(s,t)$, it is easy to see that the number of symmetries with center $x$ cannot exceed $t$, and if that number {\em is} $t$, $x$ must
be regular point. In this case, we speak of a {\em center of symmetry}. Dually, one introduces symmetries with an {\em axis}, and {\em axes of symmetry}.

\subsubsection{Moufang properties}

If $\mM$ is a subgeometry of $\Gamma$, by $\Aut(\Gamma)^{[\mM]}$ we denote the subgroup of the automorphism group $\Aut(\Gamma)$ of $\Gamma$ which fixes every line incident with a point of $\mM$ and every point incident  with a line of $\mM$. Now a (dual) root $\gamma$ is {\em Moufang} if $\Aut(\Gamma)^{[\dot{\gamma}]}$ acts transitively on the apartments containing 
$\gamma$. In fact,  $\Aut(\Gamma)^{[\dot{\gamma}]} =: A(\dot{\gamma})$ then acts {\em sharply} transitively on these apartments.
Once a (dual) root $\gamma$ is Moufang, all  (dual) roots with interior $\dot{\gamma}$ are also Moufang, with respect to the same group $A(\dot{\gamma})$. (The latter groups are uniquely defined by $\dot{\gamma}$ and the Moufang property.)
In a natural way, we also
use the terms ``Moufang i-root'' and ``dual Moufang i-root'', and the elements of $A(\dot{\gamma})$ are called {\em root-elations} throughout.

Now $\Gamma$ is {\em half Moufang} if all roots or all dual roots are Moufang. It is {\em Moufang} if all roots {\em and} dual roots are.

\subsection{Subquadrangles}

\subsubsection{Subquadrangles}

A {\em subquadrangle}\index{subquadrangle}, or also {\em subGQ}\index{subGQ},
$\mS' = (\mP',\mB',\I')$ of a GQ $\mS = (\mP,\mB,\I)$ is a GQ
for which $\mP' \subseteq  \mP$, $\mB' \subseteq  \mB$,
and where $\I'$ is the restriction of $\I$ to $(\mP' \times \mB') \cup (\mB' \times \mP')$.
The subGQ $\mS'$ is {\em ideal} if for any point $x \in \mP'$ we have that $\{ L \I x \vert L \in \mB'\} = \{ L \I x \vert L \in \mB\}$.
Dually, one speaks of {\em full subGQs}. 
\\

The following results will sometimes be used without further reference.

\bt[\cite{PTsec}, 2.2.1]
\label{2.2.1}
Let $\mS'$ be a proper subquadrangle of order $(s',t')$ of the GQ $\mS$ of order $(s,t)$. Then either $s = s'$ or
$s \geq s't'$. If $s = s'$, then each external point of $\mS'$ is collinear with the $st' + 1$ points of an ovoid
of $\mS'$; if $s = s't'$, then each external point of $\mS'$ is collinear with exactly $1 + s'$ points of
$\mS'$.
\et

\bt[\cite{PTsec}, 2.2.2]
\label{2.2.2}
Let $\mS'$ be a proper subquadrangle of the GQ $\mS$, where $\mS$ has order $(s,t)$ and $\mS'$ has order $(s,t')$
(so $t > t'$). Then we have

\begin{enumerate}
\item[$(1)$]
$t \geq s$; if $s = t$, then $t' = 1$.
\item[$(2)$]
If $s > 1$, then $t' \leq s$; if $t' = s \geq 2$, then $t = s^2$.
\item[$(3)$]
If $s = 1$, then $1 \leq t' < t$ is the only restriction on $t'$.
\item[$(4)$]
If $s > 1$ and $t' > 1$, then $\sqrt{s} \leq t' \leq s$ and $s^{3/2} \leq t \leq s^2$.
\item[$(5)$]
If $t = s^{3/2} > 1$ and $t' > 1$, then $t' = \sqrt{s}$.
\item[$(6)$]
Let $\mS'$ have a proper subquadrangle $\mS''$ of order $(s,t'')$, $s > 1$. Then $t'' = 1$, $t' = s$
and $t = s^2$.
\end{enumerate}
\et

\bt[\cite{PTsec}, 2.3.1]
\label{2.3.1}
Let $\mS' = (\mP',\mB',\I')$ be a substructure of the GQ $\mS$ of order $(s,t)$ so that the following two conditions are
satisfied:

\begin{itemize}
\item[{\rm(i)}]
if $x,y \in \mP'$ are distinct points of $\mS'$ and $L$ is a line of $\mS$ such that $x \I L \I y$, then $L \in \mB'$;
\item[{\rm(ii)}]
each element of $\mB'$ is incident with $s + 1$ elements of $\mP'$.
\end{itemize}

Then there are four possibilities:

\begin{enumerate}
\item[$(1)$]
$\mS'$ is a dual grid, so $s = 1$;
\item[$(2)$]
the elements of $\mB'$ are lines which are incident with a distinguished point of $\mP$, and $\mP'$ consists of those points
of $\mP$ which are incident with these lines;
\item[$(3)$]
$\mB' = \emptyset$ and $\mP'$ is a set of pairwise noncollinear points of $\mP$;
\item[$(4)$]
$\mS'$ is a subquadrangle of order $(s,t')$.
\end{enumerate}
\et

The following result is now easy to prove.

\bt[\cite{PTsec}, 2.4.1]
\label{2.4.1}
Let $\theta$ be an automorphism of the GQ $\mS = (\mP,\mB,\I)$ of order $(s,t)$.
The substructure $\mS_{\theta} = (\mP_{\theta},\mB_{\theta},\I_{\theta})$\index{$\mS_{\theta} =
(\mP_{\theta},\mB_{\theta},\I_{\theta})$}
of $\mS$ which consists of the fixed
elements of $\theta$ must be given by (at least) one of the following:

\begin{itemize}
\item[{\rm(i)}]
$\mB_{\theta} = \emptyset$ and $\mP_{\theta}$ is a set of pairwise noncollinear points;
\item[{\rm(i)$'$}]
$\mP_{\theta} = \emptyset$ and $\mB_{\theta}$ is a set of pairwise nonconcurrent lines;
\item[{\rm(ii)}]
$\mP_{\theta}$ contains a point $x$ so that $y \sim x$ for each $y \in \mP_{\theta}$, and each line of $\mB_{\theta}$
is incident with $x$;
\item[{\rm(ii)$'$}]
$\mB_{\theta}$ contains a line $L$ so that $M \sim L$ for each $M \in \mB_{\theta}$, and each point of $\mP_{\theta}$
is incident with $L$;
\item[{\rm(iii)}]
$\mS_{\theta}$ is a grid;
\item[{\rm(iii)$'$}]
$\mS_{\theta}$ is a dual grid;
\item[{\rm(iv)}]
$\mS_{\theta}$ is a subGQ of $\mS$ of order $(s',t')$, $s',t' \geq 2$.
\end{itemize}
\et

Finally, we recall a result on fixed elements structures of whorls.

\bt[\cite{PTsec}, 8.1.1]
\label{8.1.1}
Let $\theta$ be a nontrivial whorl about $p$
of the GQ $\mS = (\mP,\mB,\I)$ of order $(s,t)$, $s \ne 1 \ne t$.
Then one of the following must hold for the fixed element structure
$\mS_{\theta} = (\mP_{\theta},\mB_{\theta},\I_{\theta})$\index{$\mS_{\theta} = (\mP_{\theta},\mB_{\theta},\I_{\theta})$}.

\begin{enumerate}
\item[$(1)$]
$y^{\theta} \ne y$ for each $y \in \mP \setminus p^{\perp}$.
\item[$(2)$]
There is a point $y$, $y \not\sim p$, for which $y^{\theta} = y$. Put $V = \{ p,y \}^{\perp}$ and $U = V^{\perp}$.
Then $V \cup \{p,y\} \subseteq  \mP_{\theta} \subseteq  V \cup U$, and $L \in \mB_{\theta}$ if and only if $L$
joins a point of $V$ with a point of $U \cap \mP_{\theta}$.
\item[$(3)$]
$\mS_{\theta}$ is a subGQ of order $(s',t)$, where $2 \leq s' \leq s/t \leq t$, and hence $t < s$.
\end{enumerate}
\et

\medskip
\subsection{The dual net $\Pi(x)$}

The following can be found in  \cite[1.3.1]{PT}.
Let $x$ be a regular point of a thick GQ $\mathcal{S} = (\mP,\mB,\I)$ of order $(s,t)$.
Then the incidence structure $\Pi(x)$ with 

\begin{itemize}
\item
\textsc{point set} $x^{\perp} \setminus \{ x \}$;
\item
\textsc{line set} the set of spans $\{ q,r \}^{\perp \perp}$, where $q$ and $r$ are
noncollinear points of $x^{\perp} \setminus \{ x \}$, 
\end{itemize}
and with the natural
incidence, is the dual of a net of order $s$ and degree $t + 1$.
If in particular $s = t$, there arises a dual affine plane of order $s$.
Also, in the case $s = t$, the incidence structure $\widehat{\Pi(x)}$ with point set
$x^{\perp}$, with line set the set of spans $\{ q,r \}^{\perp \perp}$, where
$q$ and $r$ are different points in $x^{\perp}$, and with the natural incidence,
is a projective plane of order $s$.

\bt[\cite{notenet}]
\label{netten}
Suppose $\mathcal{S} = (\mP,\mB,\I)$ is a GQ of order $(s,t)$, $s,t \ne 1$, with a regular point $x$.
Let $\mathcal{N}_x = \Pi(x)^D$ be the net which arises from $x$, and suppose $\mathcal{N}'_x$
is a subnet of the same degree as $\mathcal{N}_x$.
Then we have the following possibilities:
\begin{enumerate}
\item[$(1)$]
$\mathcal{N}'_x$ coincides with $\mathcal{N}_x$;
\item[$(2)$]
$\mathcal{N}'_x$ is an affine plane of order $t$ and
$s = t^2$; also, from $\mathcal{N}'_x$ there arises a proper subquadrangle
of $\mathcal{S}$ of order $t$ having $x$ as a regular point.
\end{enumerate}
If, conversely, $\mathcal{S}$ has a proper subquadrangle containing the point
$x$ and of order $(s',t)$ with $s' \ne 1$, then it is of order $t$, and hence $s = t^2$. Also, there arises a proper subnet of
$\mathcal{N}_x$ which is an affine plane of order $t$.
\et

\medskip
\subsection{Main Parameter Conjecture}

The {\em Main Parameter Conjecture} (MPC) for finite generalized quadrangles (of order $(s,t)$, $s \leq t$) states that 
either $t = s$ ($t$ a prime power), $t = s + 2$ ($t - 1$ a prime power), $t = \sqrt[2]{s^2}$ ($t$ a prime power) or $t = s^2$ ($t$ a prime power).  
It is the analogue for finite generalized quadrangles of the Prime Power Conjecture for finite projective planes,  but regarded to be even more 
out of reach. 
Many famous conjectures (such as Kantor's conjecture for EGQs or Payne's conjecture for STGQs) are special cases (even 
concentrating only on the prime power part, and forgetting about the precise form of $s$ and $t$). 

There also exists a theory for the infinite case, and we refer to \cite{Order}, and the references therein, for much more information on both 
the finite and infinite case.\\

We aim at establishing MPC for the entire class of STGQs in the present series of papers.

\newpage
\vspace*{6cm}
\begin{center}
\item
{\bf PART 0}
\item
\item
{\bf SKEW TRANSLATION QUADRANGLES}
\end{center}

\newpage
\section{STGQs}

For a generalized quadrangle $\mS = (\mP,\mB,\I)$, we call a point $x \in \mS$  an {\em elation point}\index{elation!point} if there exists a group $K \leq \Aut(\mS)$ which fixes $x$ linewise and acts sharply transitively on $\mP \setminus x^{\perp}$. We call the latter set the {\em affine points} w.r.t. $\mS^x$. Similarly, the {\em affine lines} of $\mS^x$ are the lines which are not incident with $x$.
The group $K$ is the {\em elation group}\index{elation!group} of the {\em elation generalized quadrangle} (EGQ) $(\mS^x,K)$.  EGQs are the natural equivalents for generalized quadrangles as translation planes are for projective planes, although several basic questions for EGQs turn out to have different answers 
than in the theory for planes. For instance, elation groups are not necessarily abelian, and EGQs can have different (even non-isomorphic) elation groups for the same elation point.

Let $(\mS^x,K)$ be a finite thick EGQ of order $(s,t)$. Let $z \in \mP$ be not collinear with $x$, and let $\{U_i \vert i \in I\}$ be the lines incident with $z$, where 
for the index set $I$ we have $\vert I \vert = t + 1$. For each $i \in I$, let $u_i := \proj_{U_i}x$, and put $[U_i] := xu_i$. Define $\mF := 
\{ K_{U_i} \vert i \in I \}$ and $\mF^* := \{ K_{u_i} \vert i \in I \}$, and note that choosing a different point for $z$, say $z^\iota$ with $\iota \in K$, just amounts 
to taking the image of each element of $\mF \cup \mF^*$ under $\iota$. So often we do not specify the point $z$ at all. 
For $A = K_{U_j} \in \mF$, we denote $K_{u_j}$ by $A^*$ for the sake of convenience.
Note the following properties:

\begin{itemize}
\item[(a)]
$\vert K \vert = s^2t$ and for all $A \in \mF$ we have $\vert A \vert = s$ and $\vert A^* \vert = st$; also,
$\vert \mF \vert = \vert \mF^* \vert = \vert I\vert$;
\item[(b)]
for all $A \in \mF$ we have $A \leq A^*$;
\item[(c)]
for all two by two different $A, B, C$ in $\mF$ we have $AB \cap C = \{\id\}$;
\item[(d)]
for different $A, B$ in $\mF$ we have $A \cap B^* = \{\id\}$.
\end{itemize}

If an abstract group $K$ has the properties (a) through (d), we call $(\mF,\mF^*)$ a {\em Kantor family} of type $(s,t)$ in $K$. 
One then constructs an EGQ as follows: points are of three kinds | a symbol $(\infty)$ (which will be the elation point), cosets $kA^*$ for $k \in K$ and $A \in \mF$ (the points collinear with and different from $(\infty)$), and the elements of $K$ (the affine points). Lines are of two kinds | symbols $[B]$ for $B \in \mF$ (the lines on $(\infty)$) and cosets $rB$ for $r \in K$ and $B \in \mF$ (the affine lines).
Incidence is indicated by the descriptions we wrote between brackets. It is easy to show that this incidence structure $\Gamma(\mF,\mF^*)$ is an EGQ
of order $(s,t)$, for which $(\infty)$ is the elation point and $K$ serves as elation group, acting by left multiplication. Moreover, starting from an EGQ $(\mS^x,K)$ as above, one easily shows that $(\mS^x,K)$ is isomorphic to $(\Gamma(\mF,\mF^*)^{(\infty)},K)$.

\br{\rm
Only slight adjustments are needed to incorporate the infinite case.
}
\er

\subsection{Skew translation generalized quadrangles}

A {\em skew translation generalized quadrangle} (STGQ) $(\mS^x,K)$ is an EGQ with the additional local Moufang property that $K$ contains a subgroup
$\mathbb{S}$ which consists of $t$ symmetries with center $x$ (which is the theoretical maximum). This property expresses precisely the fact that 
$K_{u,v}$ with $u \sim x \sim v \not\sim u$ is independent of the choice of $u, v$ (all these groups coincide with $\mathbb{S}$).
One can show that a finite EGQ $(\mS^x,K)$ is an STGQ if and only if $x$ is a regular point. In terms of Kantor families, an EGQ is an STGQ if and only if 
for each $A, B \ne A$ in $\mF$, $A^* \cap B^*$ is independent of the choice of $(A, B) \in \mF \times \mF \setminus \mathrm{diag}$.
Another way to phrase this is asking that there is a normal subgroup $\mathbb{S}$ of $K$ such that for each $A \in \mF$, $A^* = A\mathbb{S}$; 
for each $A, B$ $(\ne A)$ $\in \mF$ we then have that $A^* \cap B^* = \mathbb{S}$.\\

Up to a combination of so-called ``Payne integration'' and point-line duality, all known finite generalized quadrangles except the classical examples $\hH(4,q^2)$ related to $4$-dimensional projective Hermitian varieties and their duals  are STGQs, and can hence be  constructed from
Kantor families as described above. 

No classification results are known for STGQs with general parameters; the only ``general'' result is the following, 
obtained independently by Chen \cite{XC} and Hachenberger \cite{Hach}.

\bt[\cite{XC,Hach}]
%\label{par}
The parameters $s,t$ of a finite STGQ $(\mS^x,K)$ are powers of the same prime.
\et

(Chen's proof was never published by him, but can be found in \cite{LEGQ}, along with many other results concerning EGQs and STGQs.)
So $K$ is a $p$-group for some prime $p$. We often denote the  prime $p$ by $\mathrm{char}(\mS^x)$.\\

Before proceeding with describing the known STGQ constructions, we introduce {\em general Heisenberg groups}.

\subsection{The general Heisenberg group}

The {\em general Heisenberg group}\index{general Heisenberg group} $\hH_n(q)$ of dimension $2n + 1$\index{general Heisenberg group!of dimension $2n + 1$} over $\mathbb{F}_q$, with $n$ a natural number,
is the group of square $(n + 2)\times(n + 2)$-matrices with entries in $\mathbb{F}_q$, of the following form (and with the usual matrix multiplication):
\begin{equation} \left(
 \begin{array}{ccc}
 1 & \alpha & c\\
 0 & \id_n & \beta^T\\
 0 & 0 & 1\\
 \end{array}
 \right),                                         \end{equation}

 \noindent
 where $\alpha, \beta \in \mathbb{F}_q^n$, $c \in \mathbb{F}_q$ and with $\id_n$\index{$\id_n$} being the $n\times n$-unit matrix.
 Let $\alpha, \alpha',\beta, \beta' \in \mathbb{F}_q^n$ and $c,c' \in \mathbb{F}_q$; then
 \begin{equation}       \left(
 \begin{array}{ccc}
 1 & \alpha & c\\
 0 & \id_n & \beta^T\\
 0 & 0 & 1\\
 \end{array}
 \right)
 \times
  \left(
 \begin{array}{ccc}
 1 & \alpha' & c'\\
 0 & \id_n & {\beta'}^T\\
 0 & 0 & 1\\
 \end{array}
 \right)
 =
  \left(
 \begin{array}{ccc}
 1 & \alpha + \alpha'& c + c' + \langle \alpha,\beta' \rangle\\
 0 & \id_n & \beta + \beta'\\
 0 & 0 & 1
 \end{array}
 \right). \end{equation}

Here $\langle x,y\rangle$\index{$\langle x,y \rangle$}, with $x = (x_1,x_2,\ldots,x_n)$ and $y = (y_1,y_2,\ldots,y_n)$ elements of $\mathbb{F}_q^n$, denotes $x_1y_1 + x_2y_2 + \cdots + x_ny_n = xy^T$.
Note that
$\hH_n(q)$ is isomorphic to the group $\{(\alpha,c,\beta) \vert \alpha,\beta \in \mathbb{F}_q^n, c \in \mathbb{F}_q\}$, where the group operation $\circ$ is given by
\begin{equation}   (\alpha,c,\beta)\circ(\alpha',c',\beta') = (\alpha + \alpha', c + c' + \alpha{\beta'}^T, c + c').                 \end{equation}

%Throughout these notes, we keep using the latter representation for the general Heisenberg group.

\bt
\label{GHG}
The following properties hold for $\hH_n(q)$ (defined over $\mathbb{F}_q$).
\begin{itemize}
\item[{\rm (i)}]
$\hH_n(q)$ has exponent $p$ if $q = p^h$ with $p$ an odd prime; it has exponent $4$ if $q$ is even.
\item[{\rm (ii)}]
The center of $\hH_n(q)$ is given by
$Z(\hH_n(q)) = \{(0,c,0) \vert c \in \mathbb{F}_q\}$.
\item[{\rm (iii)}]
For each $(\alpha,c,\beta) \in \hH_n(q)$, we have $(\alpha,c,\beta)^{-1} = (-\alpha,-c + \alpha\beta^T,-\beta)$.
\item[{\rm (iv)}]
$[\hH_n(q),\hH_n(q)] = Z(\hH_n(q))$ and $\hH_n(q)$ is nilpotent of class $2$.
\end{itemize}
\et

Define the following map
\begin{equation}
\chi: V \times V \longrightarrow Z(\hH_n(q)): (aZ(\hH_n(q)),bZ(\hH_n(q))) \longrightarrow [a,b],
\end{equation}
where $V := \hH_n(q)/Z(\hH_n(q))$. Then $\chi$ is a BAN-form of $V$, which can be seen as a vector space over $\mathbb{F}_q$.
If $n = 1$, the polar space of $\chi$ (which is the natural incidence geometry of the absolute subspaces of $\chi$) in the projective space $\mathbf{P}(V)$ coming from $V$ is a projective line $\PG(1,q)$; if $n = 2$, 
a classical quadrangle $\hW(q)$ arises.

\subsection{Examples}
\label{EX}

In this subsection, $\PG(n,q)$ | where $n$ is in $\mathbb{N} \cup \{-1\}$ and $q = p^h$ is a prime power |
is the  $n$-dimensional projective space over the field $\mathbb{F}_q$ (as a building). Its automorphism group is 
$\mathbf{P\Gamma L}_{n + 1}(q) \cong \PGL_{n + 1}(q) \rtimes \mbox{Gal}(\mathbb{F}_q/\mathbb{F}_p)$.\\

We start with describing the known examples of STGQs, and make some side remarks on the way.

\subsubsection{$\hW(q)$, $q$ odd}

The points of $\PG(3,q)$ together
with the totally isotropic lines with respect to a symplectic polarity form a GQ $\hW(q)$\index{$\hW(q)$}
of order $q$. A symplectic polarity $\Theta$ of $\PG(3,q)$ has the following canonical form:
\begin{equation} X_0Y_3 + X_1Y_2 - X_2Y_1 - X_3Y_0.\end{equation}

For each point $x$, $\hW(q)^x$ is an STGQ w.r.t precisely one elation group, which is isomorphic to the $3$-dimensional general
Heisenberg group $\hH_1(q)$ over $\F_q$.

No other STGQs of odd order $s$ are known.

\subsubsection{TGQs of even order $q$}

Each EGQ $(\mS^x,K)$ of even square order $(q,q)$ for which the elation group is abelian (in which case one calls the EGQ also a {\em translation generalized quadrangle}) is also an STGQ. (In that case, $K$ is elementary abelian.) Many such examples are known, cf. \cite{TGQ}.

All known even square order STGQs are in fact TGQs.

\subsubsection{Hermitian quadrangles $\hH(3,q^2)$}
\label{Herm}

Next, let $\hH$ be a nonsingular Hermitian variety in
$\PG(3,q^2)$. The $\F_q$-rational points and lines of $\hH$ form a generalized quadrangle
$\hH(3,q^2)$\index{$\hH(3,q^2)$}, which has order $(q^2,q)$.
The variety $\hH$ has the following canonical form:
\begin{equation} X_0 ^{q + 1} + X_1 ^{q + 1} + X_2^{q + 1} + X_3^{q + 1}   = 0.\end{equation}

For each point $x$, $\hH(3,q^2)^x$ is an STGQ w.r.t. an elation group which is isomorphic to the $5$-dimensional general
Heisenberg group $\hH_2(q)$ over $\F_q$. If $q$ is odd, there is precisely one elation group per point; if $q$ is even,
there exists a second isomorphism class of elation groups, as observed by Rostermundt \cite{Rost} and, independently, the author \cite{Basic}.

\subsubsection{Generalization: $\hH(3,q^2)$ $\mapsto$ flock quadrangles}

A {\em flock}\index{flock} $\mF$ of the quadratic cone $\mK$ in $\PG(3,q)$ is a partition of $\mK$ without its vertex into $q$ irreducible conics. 
Choosing coordinates in such a way that the equation for $\mK$ becomes $X_0X_1 = X_2^2$, one can show that the equations of the conic 
planes of the flock define a certain Kantor family in $\hH_2(q)$ of type $(q^2,q)$, and the associated quadrangle is a {\em flock quadrangle} \cite{flock}.
It is often denoted as  $\mathcal{S}(\mathcal{F})$. If one chooses a special flock for which all conic planes share some exterior line to the flock,
one can show that $\mS(\mF) \cong \hH(3,q^2)$. 

All flock quadrangles are STGQs with elation group $\hH_2(q)$, and this class of GQs was the main motivation for Payne to consider 
and study the MSTGQs in his celebrated paper \cite{9} (cf. the next section).

\subsubsection{Dual Suzuki-Tits quadrangles}

This class and its Kantor family is described in detail in the next section.

\subsubsection{Translation duals of suitable GQs}

Let $n, m$ be nonzero positive integers.
An {\em egg} $\mO = \mO(n,m,q)$ is a certain set of $(n - 1)$-spaces in $\PG(2n + m - 1,q)$ having tangent spaces at each of its elements, 
as such generalizing ovoids of $\PG(3,q)$ in a natural way, see \cite{TGQ}. Similarly as for ovoids, one can construct TGQs $\hT(\mO)$ from eggs and 
conversely, it can be shown that each TGQ is isomorphic to a $\hT(\mO)$. This natural generalization of Tits quadrangles has led to examples of 
quadrangles which weren't covered by the Tits construction. (We refer to \cite{TGQ} for all the details.)\\

Now let $n \ne m$.
The tangent spaces of an egg $\mO(n,m,q)$ in $\PG(2n+m-1,q)$ form an egg
$\mO^\ast(n,m,q)$\index{$\mO^*(n,m,q)$} in the dual space of $\PG$$(2n+m-1,q)$. So in addition to the  TGQ
$\hT(\mO)$, a TGQ $\hT(\mO^*)$\index{$\hT^*(n,m,q)$} arises which is not necessarily isomorphic to $\hT(\mO)$. The egg
$\mO^\ast(n,m,q)=\mO^\ast$\index{$\mO^*$} will be called the {\em translation dual}\index{translation!dual}
of $\mO(n,m,q)=\mO$, and $\hT(\mO^\ast)$\index{$\hT^*(\mO)$}\index{$\hT(\mO^*)$}
will be called the {\em translation dual}\index{translation!dual} of $\hT(\mO)$.\\

The most studied class of eggs in odd characteristic is the class of eggs which are ``good'' at some element, a property which is satisfied for each point of any
ovoid, and which expresses the fact that the corresponding ``good TGQ'' has many full subGQs on the line corresponding to the good element.
In that case, $m = 2n$, and moreover, the point-line dual of $\hT(\mO^*)$ is isomorphic to a flock quadrangle.\\

\begin{center}
good TGQ $\hT(\mO)$ $\overset{*}{\longrightarrow}$ translation dual $\hT(\mO^*)$ $\overset{D}{\longrightarrow}$ flock quadrangle\\
\end{center}

\medskip
In \cite{QKan}, solving a question of Kantor affirmatively, and using the machinery of \cite{TGQflock,two,SFGQ},
the author showed that this construction produced a new class of STGQs, namely the point-line duals $\hT(\mO)^D$ of $\hT(\mO)$.
Moreover, the author later showed that if $\hT(\mO) \not\cong \hH(3,q^2)$, the corresponding (unique) elation groups were new, not isomorphic 
to $\hH_2(q)$, as a direct corollary of the next theorem.

\bt[\cite{Isomflock}]
Let $\mS$ be an EGQ of order $(q^2,q)$, $q$ any prime power, with elation group $\mathcal{H}_2(q)$. Then $\mS$ is a flock quadrangle.
\et

\medskip
\subsection{Central STGQs}

Each of the examples described above has the property that the group of symmetries $\mathbb{S}$ with center the elation point of the STGQ is a central subgroup of the elation group. In general, we call STGQs with this property {\em central STGQs}. As we will see in the next section, the property of being central is equivalent to the strong local Moufang property ``(MSTGQ1)$^b$'', and allows one to control a great number of properties of STGQs. The latter claim will become clear throughout the present series of papers.

%\subsection{Special and extra-special $p$-groups}

%A $p$-group $P$  is {\em special}\index{special $p$-group} if either $[P,P] = Z(P) = \phi(P)$ is elementary abelian or $P$ itself is. Note that $P/[P,P]$ is %elementary abelian
%in that case. So we have the exact sequence

%\begin{equation}   \he \mapsto [P,P] \mapsto P \mapsto V(n,p) \mapsto \he,            \end{equation}

%where $V(n,p)$\index{$V(n,p)$} is the $n$-dimensional vector space over $\mathbb{F}_p$ and $\vert P\vert = p^n\vert [P,P]\vert$.

%If $\vert Z(P) = [P,P] = \phi(P)\vert = p$, $P$ is called {\em extra-special}\index{extra-special}.\\

%\medskip
%{\bf Example}.\quad
%The general Heisenberg group $\hH_1$ over $\mathbb{F}_p$  is extra-special.\\

\newpage
\vspace*{6cm}
\begin{center}
\item
{\bf PART I}
\item
\item
{\bf STRUCTURE THEORY}
\end{center}
\addcontentsline{toc}{chapter}{Part I - Structure theory}

\newpage

\section{Local Moufang conditions}

Before proceeding, we need the condition ``$(M)_x$''\index{$(M)_x$}. Let $\mS$ be a GQ, and let $x$ be a point of $\mS$.
Then $\mS$ satisfies $(M)_x$ if for any i-root $(x,L,y)$ ($x \I L \I y \ne x$), and any line $M \I y$, $M \ne L$,
the automorphism group of $\mS$ that fixes $x$ and $y$ linewise, and $L$ pointwise, $\Aut(\mS)_{[\{ x,L,y \}]}$,
acts transitively on  the points of $M \setminus \{y\}$. In that case, the action is sharply transitive, and the definition
is independent of $M$ (or $y$).  So $(M)_x$ holds if every i-root on $x$ is Moufang (by definition).

\subsection{MSTGQs}

Now let $(\mS^x,L)$ be a finite EGQ  with the following additionial properties:

\begin{itemize}
\item[(MSTGQ1)]
$\mS^x$ satisfies property $(M)_x$;
\item[(MSTGQ2)]
for each point $y \sim x \ne y$, $L_{[\{x,xy,y\}]} =: L(y)$ has the property that if $\alpha \in L(y)^\times$ fixes some affine line 
$U$ (which must be concurrent with $xy$ as $\alpha \ne \id$), it also fixes $U \cap xy$ linewise; 
\item[(MSTGQ3)]
no line $U \I x$ is the unique center of a triad $\{ V,W,X \}$, where $V \I x$.
\end{itemize}

Payne introduced such EGQs in \cite{9}, and calls them {\em Moufang skew translation generalized quadrangles} (MSTGQs). He shows in {\em loc. cit.}
that the point $x$ is regular, making it a center of symmetry for a subgroup $\mathbb{S} \leq L$, and $\mathbb{S} \leq Z(L)$ is 
elementary abelian. Moreover, for any $u \sim x \ne u$, $L(u)$ is also elementary abelian, and $L$ is a $p$-group. 

So $\mS^x$  is an STGQ, but it is of a very particular type: not all STGQs satisfy (MSTGQ1) (as we will see further in the next section), 
and no known square STGQ of even order has (MSTGQ3).  A more general version of (MSTGQ1) seems to be true (now supposing that $(\mS^x,L)$ is an STGQ), namely:

\begin{itemize}
\item[(MSTGQ1)$^b$]
Each $\alpha \in L^\times$ which fixes some affine line 
$U$, also fixes $U \cap \proj_Ux$ linewise.
\end{itemize}

Already from (MSTGQ1)$^b$, one can prove that the symmetry group $\mathbb{S}$ is in the center of $L$, and conversely, $\mathbb{S} \leq Z(L)$
implies (MSTGQ1)$^b$. So this seems the quintessential Moufang property for STGQs to study and obtain.

\subsection{Property (G)}

In this subsection we have a deeper look into the so-called ``Property (G)''.\\

Let $x$ be a point of a GQ $\mS$ of order $(t^2,t)$, and let $U, V$ be distinct lines incident with $x$. Then $\mS$ satisfies
{\em Property (G)} at the pair $\{U,V\}$ if any triad of lines $\{V,W,Z\}$ in $U^{\perp}$ is $3$-regular (which means that $\vert \{V,W,Z\}^{\perp} \vert =  \vert \{V,W,Z\}^{\perp\perp} \vert = t + 1$. (Note that the definition is symmetric in $U$ and $V$.)
The flag $(x,L)$ has {\em Property (G)} if all pairs of distinct lines $\{L,M\}$ ``on $x$'' have (G).
One says that $x$ has {\em Property (G)} if all pairs $\{U,V\}$ ``on $x$'' have Property (G).\\

The following theorem was first obtained in odd characteristic in \cite{III}, answering a fundamental conjecture of Payne's essay \cite{9}.
In the case of even characteristic, it was obtained only much later by M. R. Brown \cite{broeven}. (In odd characteristic, only 
one flag was required in \cite{III}.) 

\bt[\cite{broeven,III}]
\label{Geven}
A GQ of order $(t^2,t)$ satisfying Property (G) at two distinct flags $(u,L)$ and $(u,M)$ for a point $u$ is isomorphic to a flock GQ.
\et

\bc
\label{subTGQs}
A GQ $\mS$ of order $(t^2,t)$ with $t$ even is a flock GQ if and only if it contains a point $x$ which is in $t^3 + t^2$ distinct subTGQs of order $t$ for which $x$ is 
a translation point.
\ec

{\em Proof}.\quad
It is well known that $t^3 + t^2$ is the maximum number of subGQs of order $t$ that contain the same point in a GQ of order $(t^2,t)$. If this is the case for $x$, we have that 
for any triad of lines $\{U,V,W\}$ for which $U \I x$ and which has a center $C$ which is incident with $x$, $U, V, W, C$ and $D$ are contained in a subGQ $\mS'$ of order $t$, where $D$
is any center of $\{U,V,W\}$ different from $C$. Since $\mS'$ is a TGQ with translation point $x$, any line incident with $x$ is regular in $\mS'$, and so it easily follows that 
$x$ satisfies Property (G). Whence by Theorem \ref{Geven}, $\mS$ is a flock GQ.\\

Conversely, it is well known that a flock GQ of order $(t^2,t)$ with $t$ even has $t^3 + t^2$ subTGQs of order $t$ containing the special point, which is a common translation point for all these
subGQs.
\eop \\

Let $\mS$ be a GQ satisfying the same properties as in the statement of Corollary \ref{subTGQs}, but now with $t$ odd. Then by the same reasoning, $\mS$ is a flock GQ.
The only flock GQs in odd characteristic with ideal subGQs are isomorphic to the Kantor-Knuth flock  quadrangles \cite{TGQ}, and in that case, all ideal subGQs are isomorphic to
$\hW(t)$, contradicting the fact that we assumed the subGQs to be TGQs. So for $t$ odd, Corollary \ref{subTGQs} makes no sense.

\subsection{The Suzuki-Tits quadrangles}

Let $q = 2^e$ be an odd power of $2$ (so that $\sigma$ is a Tits endomorphism), and let $\sigma \in \mathrm{Aut}(\mathbb{F}_q)$ be such that $\sigma^2 = 2$. Now define a map

\begin{equation}
f: \mathbb{F}_q^2 \mapsto \mathbb{F}_q^2: (a,b) \mapsto a^{\sigma + 2} + ab + b^{\sigma}. 
\end{equation}

The {\em Tits ovoid} of $\PG(3,q)$ is then given by 

\begin{equation}
\mO = \{(0,1,0,0)\} \cup \{(1,f(a,b),a,b) \vert a, b \in \mathbb{F}_q\}.
\end{equation}

It is well known that this ovoid admits the natural $2$-transitive action of the Suzuki group $\mathbf{Sz}(q)$.
Embed $\PG(3,q)$ as a hyperplane in $\PG(4,q)$ by $(x,y,z,w) \longrightarrow (0,x,y,z,w)$.
Construct the TGQ $\Gamma = \hT_3(\mO)$, and call it the Suzuki-Tits quadrangle. Then $\mathrm{Aut}(\Gamma)_{(\infty)}$ (where $(\infty)$ is the special point of 
$\hT_3(\mO)$) admits the  natural $2$-transitive action of $\mathbf{Sz}(q)$ on the $q^2 + 1$ lines incident with $x$.
Also, $\hT_3(\mO)$ is classical if and only if $\mO$ is an ellipic quadric if and only if $q = 2$.

Now define  collineations $\theta(a,b,c,d,e): (u,x,y,z,w) \mapsto (u,x,y,z,w)[a,b,c,d,e]$, with $a,b,c,d,e \in \mathbb{F}_q$,  of $\PG(4,q)$ as follows:

\begin{equation} 
[a,b,c,d,e] = \left(
 \begin{array}{ccccc}
 1 & 0 & c & d & e\\
 0 & 1& f(a,b) &a & b\\
 0 & 0 & 1&0&0\\
 0 & 0 & a^{\sigma + 1} + b & 1 & a^{\sigma}\\
 0 & 0 & a & 0 & 1\\ 
 \end{array}
 \right).                                         \end{equation}

Then by \cite{StanPriv}, $G = \{ \theta(a,b,c,d,e) \vert a,b,c,d,e \in \mathbb{F}_q\}$ is a group of order $q^5$ with binary operation
\begin{equation}
\begin{aligned}
&
[a,b,c,d,e][a',b',c',d',e']\\ 
&= [a + a', b + b' + a{a'}^{\sigma}, c + c' + d({a'}^{\sigma + 1} + b') + ea', d + d', e + e' + d{a'}^{\sigma}].
\end{aligned}
\end{equation}

The group $G$ leaves $\mO$ invariant; it fixes the line $L$ of $\Gamma$  corresponding to $(0,0,1,0,0)$ $\in \mO$ pointwise, and acts sharply transitively on the lines 
not concurrent with $L$. So $\Gamma^D$ is an EGQ with elation group $G$. The center of $G$ is $Z(G) = \{ [0,0,c,0,0] \vert c \in \mathbb{F}_q\}$, and is a group 
of symmetries about $L$. So the dual of $\Gamma$ is a central STGQ of order $(t^2,t)$.\\

Now define for $t \in \mathbb{F}_q$:
\begin{equation}
A(t) = \{[a,b,tf(a,b),ta,tb] \vert a, b \in \mathbb{F}_q\}, A^*(t) = \{[a,b,c,ta,tb] \vert a, b, c \in \mathbb{F}_q\},
\end{equation}

and also put
\begin{equation}
A(\infty) = \{[0,0,0,d,e] \vert d, e \in \mathbb{F}_q\}, A^*(\infty) = \{[0,0,c,d,e] \vert c, d, e \in \mathbb{F}_q\}.
\end{equation}

Then $\hF = \{ A(t) \vert t \in \mathbb{F}_q \cup \{(\infty)\}\}$ and $\hF^* = \{ A^*(t) \vert t \in \mathbb{F}_q \cup \{(\infty)\}\}$ define the Kantor family 
of $G$ corresponding to the point $\langle (1,0,0,0,0),(0,1,0,0,0)\rangle$ of $\Gamma^D$.

\subsection{Some strange properties}

The following structure lemma is easy, but rather mysterious due to (ii), (iii) and (vi), cf. the remark after the observation.

\bo
\label{Structob}
\begin{itemize}
\item[{\rm (i)}]
$Z(G)$ is the symmetry group w.r.t. the elation point, and is elementary abelian.
\item[{\rm (ii)}]
$A^*(\infty)$ and $A(\infty)$ are elementary abelian and $A^*(\infty) \unlhd G$, so that $A^*(\infty)$ fixes $[A(\infty)]$ pointwise.
\item[{\rm (iii)}]
For $t \in \mathbb{F}_q$, $A^*(t)$ and $A(t)$ are non-abelian of exponent $4$; moreover, for $t \ne t'$, $A^*(t) \cong A^*(t')$ and $A(t) \cong A(t')$.
Also, no $A^*(t)$ is normal in $G$ and hence does not fix $[A(t)]$ pointwise.
\item[{\rm (iv)}]
$G$ is nonabelian of exponent $4$ and $G/Z(G)$ is not abelian.
\item[{\rm (v)}]
$G$ is the complete set of elations about $x$.
\item[{\rm (vi)}]
$\Gamma^D$ is not an MSTGQ.
\item[{\rm (vii)}]
$[G,G]$ is strictly contained in $A^*(\infty)$.
\end{itemize}
\eo

{\em Proof}.\quad
(i)-(ii)-(iii) are obvious. The fact that $G/Z(G)$ is not abelian follows from the fact that for any $t \in \mathbb{F}_q$, $A^*(t)/Z(G) \cong A(t)$.
Claim (v) is proved in \cite{Notes}. Claim (vi) follows from (v) and (iii). Finally, clearly $[G,G]$ fixes $[A(\infty)]$ pointwise so that it is a 
subgroup of $A^*(\infty)$. As $G/A^*(\infty)$ is isomorphic to $A(t)$ for any $t \in \mathbb{F}_q$, it follows that $A^*(\infty) \ne [G,G]$ by (iii).
\eop \\

Note that by (ii) and (iii), each i-root $(x,[A(\infty)],z)$ is Moufang, while no other i-root of $\Gamma^D$ is.\\

The Kantor family described in the previous paragraph is, by (ii)-(iii), the only known one for which not all elements are isomorphic.
Besides that, $\Gamma^D$ is the only known STGQ which is not an MSTGQ. The elation group  $G$ can not be constructed by the twisting method described in this paper | 
essentially there are two reasons; first of all, the fact that $\Gamma^D$ is not an MSTGQ (which is needed for making the twist), and secondly, if $\Gamma^D \not\cong \hH(3,t^2)$,
it cannot contain subGQs of order $t$ which are fixed elementwise by some involution of $\hT_3(\mO)$ (and such involutions are also needed).

\bo
Let $[\infty]$ be the special line of $\Gamma^D$ corresponding to $\mO$. Then every point incident with $[\infty]$ is an elation point for a unique elation group.
\eo

{\em Proof}.\quad
Follows from Observation \ref{Structob}(v) and the fact that $\mathrm{Aut}(\Gamma)^D$ is transitive on the lines through $x$.
\eop

\medskip
\subsection{SubTGQs of $\Gamma^D$}

Each nontangent plane section $\Pi \cap \mO = \mO'$ defines a $\hT_2(\mO')$ of Tits, and hence a full subTGQ of $\hT_3(\mO)$.
Moreover, any such $\mO'$ is a translation oval, so that $\hT_2(\mO')^D$ also is a subTGQ of $\Gamma^D$. Whence for any choice of $x \in \mO$, 
$\Gamma^D$ admits precisely $q^3 + q^2$ subTGQs of order $q$ that contain $x$.  
Still, $\Gamma^D$ is not a flock GQ (supposing that it is not classical). The reason why Corollary \ref{subTGQs} cannot be applied is reflected in the following theorem.

\bt
Let $\mO$ be non-classical.
\begin{itemize}
\item[{\rm (i)}]
Then $\Gamma^D$ contains precisely $q^4 + q^2$ distinct subTGQs of order $q$ which all contain $[\infty]$, and for which the (unique) translation point is a point of $[\infty]$. 
\item[{\rm (ii)}]
Each point $z \I [\infty]$ is
contained in precisely $q^3 + q^2$ subGQs of order $q$. Precisely $q^2$ of them have $z$ as (unique) translation point. For each other point $z' \I [\infty]$, $z' \ne z$,
there are precisely $q$ subGQs of order $q$ containing $\{z,z'\}$, and having $z'$ as (unique) translation point.\eop
\end{itemize}
\et

(The proof is left to the reader.)

\newpage
\section{Fix point theory and structure}

In this section, $(\mS^x,G)$ is an STGQ of order $(s,t)$ with $\mathbb{S}$ the associated group of symmetries about $x$.

\subsection{Property (*)}

We need to introduce the following crucial property:

\begin{quote}
{\rm (*)\ \ If $y \sim x \ne y$ is fixed by some element of $G$, it fixes $yx$ pointwise.}
\end{quote}

We will also say that (*) is satisfied for a subgroup $A \leq G$ if the property holds restricted to $A$.

\bt
\label{ab1}
An STGQ satisfies (*) if and only if $T = H/\mathbb{S}$ is abelian. 
If one of these equivalent conditions is satisfied, the members of $\hF$ are abelian.
\et
{\em Proof}.\quad
If $T$ is abelian, $\mS^x$ obviously satisfies (*), as $H$ acts transitively on the point sets $L \setminus \{x\}$, $L \I x$. Now let $\mS^x$ satisfy (*). Consider the associated Kantor family $(\hF,\hF^*)$, and let $A \in \hF$.
Define $G_A = \langle B \vert B \in \hF \setminus \{A\}\rangle$, and note that $G = G_A$. Note also that for $a \in A$ and $c \in C$, $C \in \hF \setminus \{A\}$, 
we have that $[a,c] \in \mathbb{S}$. So $[A,G_A] = [A,G]$ vanishes in $T$.  It now easily follows that $T$ is abelian.
As for the final part of the statement, note that for $A \in \hF$, $A \cong A\mathbb{S}/\mathbb{S}$.
\eop \\

Later, we will study STGQs for which all members of $\hF$ are abelian, but $H/\mathbb{S}$ is not.

\subsection{Fix point theory}

\bt[Fix point structure for STGQs]
\label{fix1}
Let $\mS = (\mS^x,G)$ be an STGQ of order $(s,t)$. For $\fix(e)$, $e$ being an element of $G^{\times}$, we have the following
possibilities.
\begin{itemize}
\item[{\rm (i)}]
$\fix(e) = x^{\perp}$;
\item[{\rm (ii)}]
$\fix(e) = \{x\}$;
\item[{\rm (iii)}]
$\fix(e)$ is a (``proper'') subset of some line on $x$;
\item[{\rm (iv)}]
If we are in {\rm (iii)}, and  $\alpha(e) := \vert \fix(e)\vert - 1$ is a power of $p := \mathrm{char}(\mS^x)$, then
$e$ fixes affine lines if $s$ is not even when $s = t$;
%\item[{\rm (v)}]
%If in the previous case $s = t^2$, that line is fixed pointwise;
\item[{\rm (v)}]
If $s = t$, $G$ satisfies Property (*);
\item[{\rm (vi)}]
If $s = t$, $e$ always has a fixed point different from $x$ (that is, any element of $H$ either fixes one or all lines on $x$ pointwise);
\item[{\rm (vii)}]
If $K(e)$ is the set of fixed affine lines of $e$, then $\vert K(e) \vert \equiv 0\mod{p}$. 
\end{itemize} 
\et

{\em Proof}.\quad
If we are not in case (ii) and case (iii), so that $e$ has fixed points $y, z \not\sim y$, then obviously
$e$ must be a symmetry about $x$ | i.e., $\fix(e) = x^{\perp}$. \\

%Let $\alpha(e)$ be as in (iii); let $\Delta$ be the set of fixed points of $e$ different from $x$ | suppose it is a subset of $U \I x$.\\

Consider case (iv). Suppose that $e$ does not fix affine lines, and suppose that $s$ is odd when $s = t$.
We want to apply \cite[1.9.1]{PT} (a result first obtained by Benson), which states that for any automorphism $\alpha$ of a GQ of order $(s,t)$, one has

\begin{equation}  (t + 1)\vert \fix(\alpha)\vert + g(\alpha) \equiv st + 1 \mod{s + t},                
\end{equation}
where $g(\alpha)$ is the number of points that are mapped onto a collinear but different point under $\alpha$.
Hence we have

\begin{equation}   (t + 1)\vert \fix(e)\vert + ts + (s + 1 - \vert\fix(e)\vert) \equiv st + 1 \mod{s + t},\end{equation}
or also

\begin{equation}  t\alpha(e) \equiv 0\mod{s + t},       
\end{equation}                    
contradiction. So (iv) follows from (iii).\\

(v) follows immediately from the plane structure of $\Gamma(x)$, as does (vi).\\ 

Finally, let $K(e)$ be as in (vii) (supposing that $K(e) \ne \emptyset$, as the empty case is trivial). Then $C_H(e)$ is nontrivial and it is easy to find prime elements in $C_H(e)$ (for instance in $C_{\mathbb{S}}(e)$) which act without trivial orbits on $K(e)$. So $\vert K(e) \vert \equiv 0\mod{p}$.
\eop \\

\br{\rm
Note that all these possibilities do occur.\\
}\er

%\bigskip
%Now let $s = t$, $a \ne \he$ an element of $G$ fixing some affine line $A$, $\alpha + 1$ be the number of fixed
%points on $B = \proj_xA$, and $\beta + t + 1$ be the number of fixed lines.
%Applying Benson's criterion, we obtain

%\begin{equation} (\alpha + \beta + 1)s + \alpha \equiv 0\mod{2s}.              \end{equation}

%Since $2s$ cannot divide $\alpha$, $\alpha + \beta + 1 \equiv 0 \mod{2}$, and whence $s + \alpha \equiv 0 \mod{s}$.
%So $a$ fixes $B$ pointwise.\\

%Slightly altering the initial assumption, we have

%\bl
%In a GQ of order $s$, for an elation about a point $x$ we have the following possibilities:
%\begin{itemize} 
%\item[{\rm (i)}]
%it fixes some affine line, and its projection on $x$ is fixed point by point;
%\item[{\rm (ii)}]
%it does not fix an affine line, and $\alpha \equiv 0 \mod{2s}$.
%\end{itemize}
%\eop 
%\el

%With the additional structure of STGQ, we could have concluded this lemma immediately by using the plane structure of $\Pi_x$.

Let $(\mS^x,G)$ be an STGQ of order $(s,t)$. Let $z \in Z = Z(G)$ be a central element, and suppose 
$z$ fixes some line $M \nI x$; then clearly any line meeting $\proj_xM$ is fixed by $z$.
Plugging this into Benson's identity (for points), one obtains 

\begin{equation}  (t + 1)(s + 1) + (s + 1)st \equiv st + 1 \mod{s + t}.               \end{equation} 

\noindent
On the other hand, since $\mS^x$ admits point-symmetries, we have the dual
\begin{equation}  (t + 1)st \equiv 0 \mod{s + t},               \end{equation}

so that $s/t + 1$ divides $\mathrm{gcd}(s + 1,t + 1)$ if $t \ne s$. In this (for our means for now only interesting) latter case,
we can put $s = p^{n(a + 1)}$ and $t = p^{na}$, $p$ being a prime and $n,a$ natural numbers.
As for $c \in \mathbb{N}_{0,1}$, $c^n + 1$ divides $c^m + 1$ if and only if $m \equiv 0\mod{n}$ and $m/n$ is odd,
we have a contradiction. So we have:

\bl
\label{zt}
\begin{itemize}
\item[{\rm (i)}]
If $s \ne t$ no element of $Z^{\times}$ sends an affine point to a collinear point. In particular, if the conditions of Theorem \ref{fix1}(iv) are also satisfied, elements of $Z^\times$ either fix only the point $x$, or are contained in $\mathbb{S}$.
\item[{\rm (ii)}]
If $s \ne t$, $\mS^x$ does not admit line symmetries.
\end{itemize}
\eop
\el

Now let $N \I x$, and $M \sim N$, $M \not\in x$. 
Let $g \in G_M^{\times}$, $\alpha$ the number of fixed points on $N$ besides $x$, and $f$ the number of fixed affine lines. Then Benson's identity for lines together with Lemma \ref{zt}(i) tells us that

\begin{equation}  (s + 1)(t + 1 + f) + \alpha t - f \equiv st + 1 \mod{s + t},                    \end{equation}

\noindent
so that 

\begin{equation} \alpha - f  \equiv 0 \mod{s/t + 1}.                  \end{equation}

First note that $f$ is a power of $p$ if (*) is satisfied. First let $y \I N$ be a point different from $x$ which is incident with precisely $f_y > 0$ fixed affine  lines for $\alpha$; clearly $f_y > 1$.
Let $\mathcal{L}_y$ be the set of fixed affine lines on $y$. Then any element of $\mathbb{S}$ which stabilizes $\mathcal{L}_y$ obviously centralizes $\alpha$, so that  $C_{\mathbb{S}}(\alpha)$ acts sharply transitively on $\mathcal{L}_y$. Whence the latter set's order is a power of $p$. This observation holds in general, that is, without assuming (*), so we state it as an intermediate observation:

\bo
Let $(\mS^x,G)$ be an STGQ, and let $\alpha \in G$ fix some affine line incident with $y \sim x \ne y$. Let $\mathcal{L}_y$ be the set of fixed affine lines on $y$.
Then $C_{\mathbb{S}}(\alpha)$ acts sharply transitively on $\mL_y$. In particular, if the STGQ is finite, $\vert \mL_y \vert = \vert C_{\mathbb{S}}(\alpha) \vert$.
\eop
\eo

Proceeding as before the observation,
if $U$ and $V$ are distinct fixed affine lines for $\alpha$, $U \not\sim V$, any element $\beta$ of  $G$ sending $U$ to $V$ centralizes $\alpha$, as $[\alpha,\beta]$ fixes $V$ and is contained in $\mathbb{S}$ by (*).

\medskip
\subsection{Averaging lemma}

The following lemma will be used many times.

\bl[Averaging lemma]
\label{AV}
Let $(\mS^x,H)$ be an STGQ of order $(s,t)$ with (*), where $s \ne t$ if $s$ is even. Let $L$ be an affine line. Then each element of $H_L^\times$ fixes precisely
$s$ affine lines (and so $t + s + 1$ lines in total). As a corollary, $f(\gamma)/\alpha(\gamma)$ is independent of the choice of $\gamma$ if $f(\gamma) \ne 0$. (In
the case $s = t$ odd the latter condition is redundant.)
\el

{\em Proof}.\quad
First let $s = t$ be odd.
Let $\alpha \in H_L^\times$ fix $f$  affine lines. As (*) holds, we have that 
\begin{equation}
\vert \alpha^H \vert \leq \vert [H,H]\vert \leq \vert \mathbb{S}\vert = s,
\end{equation}
so that $\vert C_H(\alpha)\vert \geq s^2$. As $\vert C_H(\alpha)\vert = \vert \mP_{\alpha}\vert$ by earlier considerations, it follows that $f \geq s$. This was already 
obtained for the other parameters under consideration.
(Note that the same reasoning can be used to prove it in an alternative fashion for these parameters.)\\

Let $\gamma \in H_L^\times$; we know that $\gamma$ fixes precisely $sp^{e(\gamma)}$ affine lines, for some natural number $e(\gamma)$, $p$ being $\mathrm{char}(\mS^x)$.
Let $X$ be the set of lines in $\proj_xL^\perp$ which are not incident with $x$. Since $H$ does not contain line symmetries, $H$ acts faithfully on $X$. Of course this 
action is transitive, so that the average number of fixed lines in $X$ of an element of  $H$ is $1$. By (*), the elements of $H$ fixing elements of $X$ precisely  are 
those of $H_u$, with $u$ an arbitrary point of $\proj_xL \setminus \{x\}$. ($H_u$ has size $st$ and fixes $\proj_xL$ pointwise.) Let $sp^{\overline{e}}$
be the average number of elements of $X$ fixed by elements of $H_u^\times$. Then we have 
\begin{equation}  
1 = ((st - t)sp^{\overline{e}} + st)\frac{1}{s^2t},        
 \end{equation}

so that $\overline{e} = 0$, and hence each element of $H_u^\times$ fixes precisely $s$ affine lines.
\eop

\br
{\rm
The Averaging Lemma will provide us much combinatorial information once we are ready to pin down the isomorphism classes of STGQs. }
\er

The following corollary is important for many purposes.

\bc
\label{centrals}
Let $(\mS^x,H)$ be as in Lemma \ref{AV}. Then $\vert C_H(\gamma)\vert = s^2$, so that $[\gamma,H] = \mathbb{S}$, for $\gamma \in \cup_{u \in x^{\perp}\setminus \{x\}}H_u$.
\ec

{\em Proof}.\quad
Let $\gamma$ be as in the statement. Then $\gamma$ fixes precisely $s$ affine lines. Call the set of these lines $\mathcal{X}$.
If $\beta \in H$ sends some line $M \in \mathcal{X}$ to $M^\beta \in \mathcal{X}$, then clearly $[\gamma,\beta] = \he$. And to obtain the latter property,
$\beta$ must stabilize $\mathcal{X}$. Besides that, (*) imposes the fact that $H_L$ is abelian for any affine $L$, and whence if $U \nI x$ is such that $\gamma \in H_U$ (such $U$ exists), $H_U \leq C_H(\gamma)$. The corollary easily follows from the fact that $H$ is transitive on $\mathcal{X}$.
\eop

\medskip
\subsection{Transfer}

The next theorem is formulated in terms of EGQs, and is quite general.

\bt[Transfer]
\label{Tran}
Let $(\mS^x,G)$ be an EGQ of order $(s,t)$, with $s > t$, and suppose that $st$ is a power of the prime $p$.
Let $a \in G^{\times}$ be such that $a$ fixes only points of some line $L_a \I x$,
and suppose that 
$\alpha(a)$, where $\alpha(a) + 1$ is the number of fixed points of $a$ on $L_a$, 
is a power of $p = \mathrm{char}(\mS^x)$.
\begin{itemize}
\item[{\rm (i)}]
Then $f(a)$, the number of  fixed affine lines, is not zero.
\item[{\rm (ii)}]
Suppose $f(a)$ is also a power of $p$. Then

\begin{equation}  
\frac{\alpha(a)}{f(a)} = (\frac{s}{t})^{o(a)},       
 \end{equation}
where  $o(e) \in 2\mathbb{Z}$. 
\end{itemize}
%It follows that also

%\begin{equation}  \frac{\alpha(a)}{r(a)}          \end{equation}
%is independent of the choice of $A$.
\et
{\em Proof}.\quad
Applying Benson's identity (for points) on $a$, we obtain
\begin{equation}
(t + 1)\fix(a)  + st + f(a)s + (s + 1 - \fix(a)) \equiv st + 1 \mod{s + t},
\end{equation}
from which we deduce that 
\begin{equation}
\alpha(a) - f(a) \equiv 0\mod{\frac{s}{t} + 1}.
\end{equation}
Pari (i) follows.\\

From now on, we suppose that $f(a)$ is also a power of $p$.\\

Let $f(a) \geq \alpha(a)$; as both are powers of $p$ by assumption, $f(a)$ divides $\alpha(a)$. We conclude that 
\begin{equation}
\frac{\alpha(a)}{f(a)} - 1 \equiv 0\mod{\frac{s}{t} + 1}.
\end{equation}

It is well known that for positive integers $c > 1, m > 1, n > 1$, we have that $\mathrm{gcd}(c^m + 1,c^n - 1)$ either is $1$, $2$ or 
$c^d + 1$, in which case $n/d$ is even. \\

As $\frac{s}{t} + 1 > 2$, with $\frac{\alpha(a)}{f(a)} = p^v$ and $\frac{s}{t} = p^u$ we get that $u = \mathrm{gcd}(v,u)$
and $v/u$ is even. So 

\begin{equation}
\frac{\alpha(a)}{f(a)} = (\frac{s}{t})^{o(a)} 
\end{equation}
with $o(a)$ even.\\

The case $\alpha(a) \leq f(a)$ is obtained by reversing the roles of $\alpha(a)$ and $f(a)$.
\eop

%Observation \ref{Tran} seems to indicate that the following conjecture might be true:\\

%\textsc{Conjecture}.\quad
%{\em Let $\mS$ be an EGQ with associated Kantor family $\mJ$. Then all members of $\mJ$ are isomorphic.}

%\bc
%If, under the conditions of Theorem \ref{Tran}, (*) is satisfied for $(\mS^x,G)$, each nonidentity element of an affine line stabilizer in $G$ fixes the same %number of lines.
%\eop \\
%\ec

\newpage
\section{EGQs of order $(t^2,t)$ with (*) are STGQs}

In this section, we prove the surprising fact that if $(\mS^x,H)$ is an EGQ with parameters $(t^2,t)$ such that (*) is satisfied for $H$,  then it is an STGQ.
We use the fact that $t$ is a prime power (consult the monograph \cite{LEGQ}, or see the Chen-Hachenberger result). \\

We first make some considerations under the assumption that only the members of $\mF$ are abelian ($(\mF,\mF^*)$ being the associated Kantor family).

Let $a$ be an element of $H^\times$ with only fixed points of the line $L_a \I x$. By (*), the number of such points is $s = t^2$. Let $U$ be any line 
not concurrent with $L_a$; then the triad $\{ L_a,U,U^a\}$ has precisely $t + 1$ centers, and all these centers are obviously fixed by $a$. Note that 
it follows that there are at most $t^2$ fixed affine lines, and that not all such lines are concurrent.
Now let $X$ and $Y$ be arbitrary nonconcurrent lines in $K(a)$, the set of fixed affine lines of $a$. Let $V$ be a line not incident with $x$
which meets both $X$ and $Y$, and let $\mu \in H_V$ map $X$ to $V$. Then $[a,\mu] = \id$, as $[a,\mu]$ fixes $L_a$ and 
$\proj_xV$ pointwise, while $U$ is fixed. 
So $\mu \in C_H(a)$, and from this observation it easily follows that $C_H(a)$ acts transitively on $K(a)$, readily implying that $f(a) = \vert K(a)\vert$ is 
prime power.
So the transfer formula for $a$ can be applied:

\begin{equation} 
\frac{s}{f(a)} = \frac{\alpha(a)}{f(a)} = (\frac{s}{t})^{o(a)} = t^{o(a)}         
 \end{equation}
with $o(a) \in 2\mathbb{Z}$. As $f(a) \leq t^2 = s$, $o(a) \geq 0$. As $f(a) > 1$, $o(a) = 0$ and so $\alpha(a) = s = f(a)$ for all nontrivial $a$ which 
fix points besides $x$. 

\br{\rm
We have the same arithmetic information at hand regarding the number of fixed lines of line stabilizers as for (*)-STGQs of order $(t^2,t)$. It might be possible to pass from this  point to the conclusion of Theorem \ref{arg2} below.}
\er 

%As such $a$ can only fix affine lines on $t$ points of $L_a$, it follows that all these points are fixed linewise.
%We have obtained (MSTGQ1)$^b$.

%Now let $b \sim x \ne b$, $c \sim c \ne c$, $b \not\sim c$, $y \in \{ b,c \}^{\perp} \setminus \{x\}$,
%and let $B := by$ and $C := cy$.  It is easy to show that $\langle H_B,H_C \rangle = H$, and obviously 
%$[H_B,H_{b,c}] = [H_C,H_{b,c}] = \{\id\}$, so $H_{b,c} \leq Z(H)$. 
%It follows that $H_{b,c}$ also satisfies (*). Now by the same argument as used in the proof of Theorem \ref{arg2}, it follows that $H_{b,c}$ is a 
%group of point symmetries of size $t$ in $H$.

%We have proved the next theorem.

%\bt[EGQs of order $(t^2,t)$ with (*) - A]
%\label{arg1}
%An EGQ $(\mS^x,H)$ of order $(t^2,t)$ satisfying (*) for each member of $\hF$, where $(\hF,\hF^*)$ is the associated Kantor family, is an STGQ. In particular, $x%$ is a regular point.\eop \\
%\et

%\br[Shortcut]{\rm
%There is alternative way to end the above argument, using a recent result of the author. Once one has obtained the local Moufang property, \cite{localhalf}
%implies that the dual i-roots containing $x$ are Moufang (the corresponding root groups being contained in $H$). Now use the argument of 
%Theorem \ref{arg2} below.
%}
%\er

\bigskip
If we now assume that the entire elation group has (*), a short proof of the following is possible.

\bt[EGQs of order $(t^2,t)$ with (*)]
\label{arg2}
An EGQ $(\mS^x,H)$ of order $(t^2,t)$ satisfying (*) is an STGQ. In particular, $x$ is a regular point.
\et
{\em Proof}.\quad
Let $b \sim x \ne b$, $c \sim x \ne c$, $b \not\sim c$, and consider $H_{b,c}$. Suppose that 
some element $\nu$ of $H_{b,c}$ is not a symmetry; then there is an affine line $U$ such that 
$U \not\sim U^{\nu}$. But then any of the $t + 1$ centers of $\{ B,U,U^{\nu} \}$ is fixed by $\nu$, forcing $\nu$ to be the identity. 
So $H_{b,c}$ only consists of symmetries, and $(\mS^x,H)$ is an STGQ.
\eop  \\

Theorem \ref{arg2} naturally generalizes to the infinite case, but such considerations will be done in a later paper.

\newpage
\vspace*{6cm}
\begin{center}
\item
{\bf PART II}
\item
\item
{\bf EXISTENCE OF ROOT-ELATIONS}
\end{center}
\addcontentsline{toc}{chapter}{Part II - Existence of root-elations}

\newpage
\section{Root-elations}

In this section we will show that (many) pseudo root-elations exist for certain STGQs.
In the entire section, we will be mainly concerned with STGQs satisfying (*). Some intermediate lemmas will be stated in such a way that they will also be valid in the non-(*) case, which will be handled near the end of the paper.  At several stages, we will obtain results in multiple ways, each  adding additional insight in 
the structure of $H$ (and $\mS^x$). The different approaches often generalize to different more global properties.
We will also obtain information about the exponents of skew elation groups as a corollary of the obtained results.

\subsection{Existence of root-elations, I}
\label{rootel}

In this paragraph we suppose (*) for the STGQ $\mS^x$. In a first instance, we let $s \ne t$ if $s$ is even, although we will consider the case of even $s = t$ in another paper of the series.\\

Let $L$ be an arbitrary affine line, and put $[L] = \proj_xL$; let $l \I L \I l \not\sim x$, and let $[M] \ne [L]$ be arbitrary on $x$. Put $M = \proj_l[M]$.
Let $\alpha$ be any element of $H_L^\times$; then (*) implies $H_L \leq C_H(\alpha)$. By Corollary \ref{centrals}, $\alpha$ fixes precisely $s$ affine lines (which are in 
$[L]^{\perp}$). Let $\mB_\alpha$ be the set of these lines, and let $\mP_\alpha$ be the set of $s^2$ points ``in'' $\mB_{\alpha}$. Note that for any $m, m' \I [M] \I m, m' \ne x$, we have 
$H_m = H_{m'}$. The following is obvious, and will be used without further notice (it is only valid because $\mS^x$ satisfies (*)! | it is also true when $s = t$ is even):

\bl
\label{palfa}
$\gamma \in H$ centralizes $\alpha$ if and only if $\gamma$ sends at least one point of $\mP_{\alpha}$ to some other point of $\mP_{\alpha}$. If this is the case, it stabilizes
$\mP_{\alpha}$ (and $\mB_{\alpha}$) globally.
\eop
\el

It is easy to see that $C_H(\alpha)$ ``splits over $H_L$ and $H_m$'' | that is, $C_H(\alpha) = (C_H(\alpha) \cap H_L)(C_H(\alpha) \cap H_m) = H_LC_{H_m}(\alpha)$
(for instance by using Lemma \ref{palfa}, and using the respective orders of the groups under the consideration in the last equality, cf. Corollary \ref{centrals}).
Let $[\alpha]$ be the subset of $[L]$ of points incident with a line of $\mB_{\alpha}$.

\bl
\label{8.8}
$C_{H_m}(\alpha)$ acts transitively on $\mB_{\alpha}$, so also on $[\alpha]$. Whence each point of $[\alpha]$ is incident with a constant number of fixed affine lines under $\alpha$.
This number is $s/\vert [\alpha]\vert =: r > 0$, and also equals $C_\mathbb{S}(\alpha)$.
\eop
\el

\begin{remark}{\rm
Lemma \ref{8.8} is also valid for the $s = t$ even case.}
\end{remark}

We will show that (*) forces $C_{H_m}(\alpha)$ to split over $\mathbb{S}$ and $H_M$.
To do this, put $U = M^{\alpha^{-1}}$, and let $k = \vert H_U \cap H_M\vert = \vert \{[L],M,U\}^{\perp}\vert$.  Note that $H_U \cap H_M$ is precisely $C_{H_M}(\alpha)$.
Since $H_M$ is abelian, $H_U \cap H_M$ fixes the $s/k$ affine lines on $\proj_{[M]}l^{\alpha^{-1}}$ which also meet $\{[L],M\}^{\perp}$. By symmetry, $H_U \cap H_M$ also fixes (at least) $s/k$ affine lines on $\proj_{[M]}l$. 
By way of contradiction, suppose that $k < \vert [\alpha]\vert$. Then $k < s/r$. So each element of $C_{H_M}(\alpha)$ fixes strictly more affine lines through ${\proj_{[M]}}l$ than $\alpha$ does through $L \cap [L]$ | call this property (L). 

Now suppose that the value $r = r(\alpha)$ is maximal for the function
\begin{equation}
r: \cup_{Y \in \mF}Y^{\times} \longrightarrow \{1,\ldots,t\};
\end{equation}
then by (L), we have that $k = s/r(\alpha) = \vert [\alpha]\vert$. Whence 

\begin{equation}   C_{H_m}(\alpha) = (C_{H_m}(\alpha) \cap \mathbb{S})(C_{H_m}(\alpha) \cap H_M) = C_{\mathbb{S}}(\alpha)C_{H_M}(\alpha).    \end{equation}

Note that $s/r$ (and so $\vert C_{H_M}(\alpha)\vert$) is independent of the choice of $M$ (and the choice of $l \I L$). This implies that for any affine line $V \I l$, we have that $\vert \{[L],V,V^{\alpha^{-1}}\}^{\perp}\vert = s/r$, 
in other words, $V$ contains precisely $s/r$ points of $\mP_{\alpha}$. An easy counting argument leads to the fact that $r(\alpha) = t$.
Let $\beta \in C_{H_M}(\alpha) = H_M \cap H_U$;   as $r(\beta)$ (obvious notation) $\geq s/k = r(\alpha) = t$, it follows that $r(\beta) = t$. 

\bp[Transfer]
Let $\nu \in A^{\times}$, $A \in \mF$. 
If $r(\nu) = t$ and $[\nu,\mu] = \id$, with $\mu \in B^{\times}$ and $B \in \mF \setminus \{A\}$, then $r(\mu) = t$.
\eop
\ep

Now suppose $\gamma \in \cup_{Y \in \mF}Y^{\times}$ is such that $r(\gamma)$ is the second largest value for $r(\cdot)$, and suppose furthermore that $r(\gamma) < t$. Suppose that $\gamma \in C \in \mF$.
Let $\omega \in D^{\times}$, $D \in \mF \setminus \{C\}$, be such that $[\gamma,\omega] = \id$. If $r(\omega) > r(\gamma)$, then $r(\omega) = t$ by assumption, and then we have seen that $r(\gamma) = t$ by Transfer, contradiction. So $r(\omega) \leq r(\gamma)$ for all such $\omega$, implying equality for all such $\omega$. But then we have seen that this forces $r(\gamma) = t$,
contradiction.

The following structural result follows.

\bt[Centrality $\equiv$ (MSTGQ1)$^b$]
Let $(\mS^x,H)$ be an STGQ of order $(s,t)$, where $s$ is odd if $s = t$. Suppose (*).
\begin{itemize}
\item[{\rm (i)}]
For any affine line $L$, $H_L = H_{[\proj_Lx]}$, that is, each element of $H_L$ fixes $\proj_Lx$ linewise. 
\item[{\rm (ii)}]
For each element $\alpha$ of $H_L^\times$, $\vert [\alpha]\vert = s/t$; in particular, if $s = t$ is odd, the affine  fixed line structure of $\alpha$ precisely consists
of the $t$ lines through $\proj_Lx$.
\item[{\rm (iii)}]
$Z(H) = \mathbb{S}$.  \eop \\
\end{itemize}
\et

\subsection{Exponents}
\label{Exp}

In this subsection, we first suppose that if $s \ne t$, (*) holds, and that the elements of $\hF$ are elementary abelian. 
Put $p = \mathrm{char}(\mS^x)$.

\begin{lemma}
\label{calc}
Let $G$ be a nilpotent group of class $2$, and let $a,b,c \in G$.
\begin{itemize}
\item[{\rm (i)}]
For any natural number $m$, we have that $[a^k,b] = [a,b]^k$.
\item[{\rm (ii)}]
For any natural number $m$, we have that

\begin{equation} (acb)^m = a^mc^mb^m[b,a]^{m(m - 1)/2}[c,a]^{m(m - 1)/2}[b,c]^{m(m - 1)/2}.\end{equation}
\end{itemize}
\end{lemma}

{\em Proof}.\quad
The first part is easy. The second part follows by (i) and induction. \eop \\

Take any element $s$ of $\mathbb{S}^\times$, and consider any point $y \sim x \ne y$. By the previous subsection, there exist $a \in H_{[y]}^\times$. Choose $a$ 
 to have order $p$, and note that $[s,a] = \he$. Under the hypothesis, we have that $sa$ fixes some affine line, so

\begin{equation}   \he = (sa)^p = a^ps^p = s^p,    \end{equation}  

so $\mathbb{S}$ has exponent $p$. When $p = 2$, it follows that $\mathbb{S}$ is elementary abelian. 

\bp
Let $(\mS^x,H)$ be an STGQ of order $(s,t)$. Let $(\hF,\hF^*)$ be the associated Kantor family.
 Suppose that the elements of $\hF$ are elementary abelian, and that either $s = t$ is odd, or (*) holds if $s \ne t$. Then the following properties hold.
 \begin{itemize}
 \item[{\rm (i)}]
 The elements of $\hF^*$ have exponent $p$. In particular, the exponent of
 $\mathbb{S}$ is $p$.
 \item[{\rm (ii)}]
 If $st$ is odd, $H$ has exponent $p$. If $st$ is even, it has exponent $4$.
 \end{itemize}
\ep
{\em Proof}.\quad
Part (i) has already been obtained, so we only deal with part (ii). 
As $H/\mathbb{S}$ is abelian and generated by elementary abelian groups (namely the groups $A\mathbb{S}/\mathbb{S}$ with $A \in \mF$, 
$(\mF,\mF^*)$ the associated Kantor family), it is also elementary abelian. (Note that this already implies that for any $\mu \in H$, 
$\mu^{p^2} = \id$, so the exponent of $H$ is either $p$ or $p^2$.)
Let $A, B \ne A$ be in $\mF$; then $H = AB\mathbb{S}$. 
Let $\nu$ be arbitrary in $H$, so that we can write $\nu = abs$ for $a \in A, b\in B, s \in \mathbb{S}$. As $[H,H] \leq \mathbb{S} \leq Z(H)$ by the 
assumptions and the previous subsection, $H$ has nilpotence class $2$.  So by Lemma \ref{calc}, 
\begin{equation}
\nu^p = (abs)^p = a^pb^p[b,a]^{p(p - 1)/2} = [b,a]^{p(p - 1)/2}.
\end{equation}
If $p$ is odd, the latter term is $\id$, so then $H$ has exponent $p$. 
\eop \\

The case $s = t$ even will be handled in a separate paper, and needs a different, indirect, approach. It is much harder.

\subsection{Existence of root-elations, II}

In this subsection we suppose that for the STGQ $\mS^x$ with Kantor family $(\mF,\mF^*)$, (*) holds. Moreover, we suppose that the elements of $\mF$
are elementary abelian. 
If $s = t$ (when (*) is automatically satisfied), we furthermore suppose that $s$ is odd | as in most parts of this section, we will handle the even square case separately. The goal is to obtain nontrivial root-elations in an entire different manner as before. (The results will be slightly weaker at the level of root-elations.) 
To that end, we will obtain some group theoretical results which have independent interest.

Before proceeding, we first recall the following result by Alperin and Glauberman, of \cite{AlpGlau}.

\bt[\cite{AlpGlau}]
\label{AlpGlau}
If $A$ is an abelian subgroup in a (finite) $p$-group $P$ of exponent $p$ and class at most $p + 1$, then there is a normal abelian subgroup $B$ of the same order contained in the normal closure of $A$ in $P$.
\et

\subsection*{Case $s > t$}

First suppose $s \ne t$, and let $A = H_L$, $L$ an affine line. By assumption, $A$ is elementary abelian, so that Theorem \ref{AlpGlau} and \S \ref{Exp} imply that $\langle A^H\rangle$
contains an elementary abelian subgroup $\mathbb{A}$ of order $s > t$, {\em which is normal in $H$.} We note the following properties:

\begin{itemize}
\item
$\mathbb{A}$ fixes $\proj_xL$ pointwise (by (*)).
\item
Since $\vert \mathbb{A}\vert > t$ and (*) holds, the previous point implies $\mathbb{A}$ contains elements $\alpha$ fixing some affine line in ${(\proj_xL)}^{\perp}$, say $M$.
\item
It follows that $\langle \alpha^H\rangle \leq \mathbb{A}$.
\end{itemize}

Let $\proj_Mx = m$.
Let $\eta$ be any nontrivial symmetry, and put $M^\eta = M'$. Then $\alpha\eta$ is an element of $H_m \setminus \mathbb{S}$ which sends $M$ to $M'$. 
As $\alpha\eta$ fixes $xm$ pointwise, it must fix some affine line $M''$ meeting $xm$, so that clearly there must be an $h \in H$ such that 
\begin{equation}   
\alpha^h = \alpha\eta.      
\end{equation}

(Note that $\alpha$ and $\alpha\eta$ have the same action on the lines incident with $x$.) So $[\alpha,h] = \eta \in \mathbb{A}$. We have shown

\bl
$\mathbb{S} \leq \mathbb{A}$. As a corollary, $\mathbb{S}$ is elementary abelian.
\eop 
\el

\bt[Existence of root-elations (for $s \ne t$, (*))]
Let $(\mS^x,H)$ be an STGQ of order $(s,t)$ for which (*) holds, and such that $s \ne t$ if $s$ is even. Suppose furthermore
that affine line stabilizers are elementary abelian. Then for any $L$ not incident with $x$,
$H_L$ contains at least $s/t > 1$ elements which fix $\proj_Lx$ linewise, that is, $H_L$ contains at least $s/t$ root-elations with i-root $(x,x\proj_Lx,\proj_Lx)$.
\et

{\em Proof}.\quad
Put $\proj_Lx = l$.
Since $\mathbb{A}$ is a subgroup of $H_l$ and contains $\mathbb{S}$, $\vert \mathbb{A}_L\vert = s/t$. The proof follows from the fact that $\mathbb{A}$ is abelian, so that 
$[\mathbb{A}_L,\mathbb{S}] = \{\id\}$. \eop \\

\subsection*{Square STGQs}

Now suppose that $s = t$ is odd, and let $L$ be an affine line.
Consider the nontrivial group $Z(H) \leq \mathbb{S}$; then $\langle Z(H),H_L\rangle$ is an elementary abelian group of order $sr > s$, with $r = \vert Z(H)\vert$. Now carry out the 
same argument as for the case $s \ne t$ with $A = H_L$, to conclude that 

\begin{itemize}
\item
$\mathbb{A}$ is an elementary abelian normal subgroup of $H$ of order $sr > s$;
\item
$\mathbb{S} \leq A$, so that $\mathbb{S}$ is abelian.  
\end{itemize}

\bt[Existence of root-elations (for $s = t$ odd)]
Let $(\mS^x,H)$ be an STGQ of order $s$ with $s$ odd.
For any $L$ not incident with $x$,
$H_L$ contains at least $\vert Z(H)\vert > 1$ elements which fix $\proj_Lx$ linewise, that is, $H_L$ contains at least $\vert Z(H)\vert$ root-elations with i-root $(x,x\proj_Lx,\proj_Lx)$.
\eop 
\et

To eventually complete the square GQs, we need to consider the case $s = t$ even, in which we cannot employ the property that elements of 
$\cup_{A \in \mF}A^* \setminus \mathbb{S}$ necessarily fix affine lines.

\subsection{Abelian subgroups of metabelian groups}

Recall that a group $G$ is called {\em metabelian} if its first derived group is abelian. So $G$ is metabelian if and only if there is a normal 
abelian subgroup $A$ such that $G/A$ also is abelian.

\bt[J. D. Gillam \cite{Gillam}]
\label{Gillam}
Let $G$ be a finite metabelian group, and let $B$ be an abelian subgroup of $G$ of maximal size. Then $G$ contains a normal abelian
subgroup $C \leq \langle B^G\rangle$ for which $\vert B \vert = \vert C\vert$.
\et

Following the (same) proof of A. Mann \cite{Mann1}, one obtains the following variation on Theorem \ref{Gillam}.

\bt
\label{metalemma}
Let $G$ be a finite metabelian group, and suppose $N$ is a normal subgroup of $G$.
Let $B$ be an abelian subgroup of $N$ of maximal size (in $N$). Then $G$ contains a normal abelian
subgroup $C \leq \langle B^G\rangle \leq N$ for which $\vert B \vert = \vert C\vert$.
\et

{\em Proof}.\quad
(Sketch.)
Let $G$ and $N$ be as in the statement of the theorem. Choose $B$ among the abelian subgroups in $\langle B^G\rangle \leq N$ with size $\vert B\vert$   
in such a way that it has maximal intersection with $[G,G]$. Assume that it is not normal in $G$ (by way of contradiction).
Then $N_G(B) \ne G$, so the normalizer condition implies that we can choose an element $h$ in $N_G(N_G(B)) \setminus N_G(B)$.
Define $A = BB^h$, and note the following properties:

\begin{itemize}
\item
$B^h \ne B$ so that $A \ne B$;
\item
$A \leq N$;
\item
$B$ and $B^h$ normalize each other, so $A = B[B,h]$ has class $2$ and $Z(A) = B \cap B^h$ ($B$ is maximal abelian in $A$).
\end{itemize}

Define $E = Z(A)(B \cap [G,G])[B,h]$; then $E \leq \langle B^G\rangle \leq N$, and $\vert B\vert \leq \vert E\vert$ (note that $A = BE$, and $B \cap B^h \leq E \cap B$). Since $E$ is abelian,
equality holds. But $E \cap [G,G] > B \cap [G,G]$ (as $[B,h] \not\in B$), contradiction.
\eop \\

Now let $(\mS^x,H)$ be an STGQ with (*), and metabelian $H$. Note that when (*) is satisfied, $H$ is metabelian if and only if
$\mathbb{S}$ is abelian. We have proved that this is already the case if $s$ is even and $s \ne t$ (assuming that affine line stabilizers are elementary abelian).

Choose any affine line $L$; let $u$ be such that $L \I u \sim x \ne u$. Then $H_u$ is normal in $H$. Let $B$ be abelian of maximal size in $H_u$;
then $\vert B \vert > s$ since $H_L$ is abelian and $Z(H) \leq \mathbb{S} \leq H_u$, while $\mathbb{S} \cap H_L = \{ \id\}$. So by Theorem \ref{metalemma}, $B$ can be taken to be 
normal in $H$. By (*), each element $b$ of $B \setminus \mathbb{S}$ fixes some affine line in $xu^{\perp}$, and so $[b,H] = \mathbb{S} \leq B$. 
On the other hand, $[b,\mathbb{S}] = \{\he\}$ implies that $b$ is a root-elation (with i-root $(u,xu,x)$). 
The following theorem follows immediately:

\bt
Suppose $(\mS^x,H)$ satisfies (*), and let $\mathbb{S}$ be abelian. Then for each point $y \sim x \ne y$, $\vert H_{[y]} \vert > s/t$. 
In particular, if $s$ is even with $s \ne t$ and affine line stabilizers are elementary abelian, we have this conclusion.
\eop
\et

\begin{remark}[Metabelian $H$ in even characteristic]
{\rm
Suppose that $\mathbb{S}$ is an elementary abelian $2$-group, and let $\vert \mathbb{S}\vert = 2^{n} = t$. We see $\mathbb{S}$ as an affine $n$-space $\Pi$ over $\mathbb{F}_2$.
Then we have a natural homomorphism
\begin{equation}  \eta: H \mapsto \mathbf{AGL}_n(2)          \end{equation}  
defined by conjugation; its kernel is $C_H(\mathbb{S})$. It might be interesting to study the action of the involutions in $H$ on the space $\mathbb{S}$.
}\end{remark}

\medskip
\subsection{STGQs of order $(t^2,t)$}

For STGQs of order $(t^2,t)$ with (*), there is a general combinatorial approach, independent of the characteristic, and inspired by the beginning of this section, which directly leads to centrality of $\mathbb{S}$.

\bt
\label{centsqsq}
An STGQ of order $(t^2,t)$ with (*) is central.
\et
{\em Proof}.\quad
Let $L \nI x$, and let $U \I x$ be arbitrary but different from $\proj_xL$. Put $Y = \{U,L\}^{\perp} \setminus \{\proj_xL\}$. Let $r$ be any point of $\proj_xL$ different
from $x$ and not incident with $x$. Let $R$ be any line on $r$ not on $x$; then $R$ meets $Y$ in $t$ different points (as the order of $\mS^x$ is $(t^2,t)$). Clearly, $H_L$ acts transitively on the lines incident with $r$ and different from $\proj_xL$. So if $\beta \in H_L$ fixes $R$, the fact that $H_L$ is abelian implies that $\beta$ fixes $r$ linewise.
In particular, $H_L \cap H_R$ yields a group of $(r,\proj_xL,l)$-elations of size $t$ ($l = x^\perp \cap L$).\\

Now let $\alpha \in H_L^\times$ be arbitrary. Let $V \in Y$ be arbitrary; then $\{\proj_xL,V,V^\alpha\}^{\perp}$ has size $t + 1$, two of its elements being $U$ and $L$.
Let $L' \in \{\proj_xL,V,V^\alpha\}^{\perp} \setminus \{U,L\}$; then $\alpha \in H_{L'} \cap H_L$, so that the first part of the theorem yields that $\alpha$ fixes $l$ linewise. 
\eop \\

\newpage
\vspace*{6cm}
\begin{center}
\item
{\bf PART III}
\item
\item
{\bf SQUARE ODD STGQS | CLASSIFICATION}
\end{center}
\addcontentsline{toc}{chapter}{Part III - Square odd STGQs | classification}

\newpage
\section{Square STGQs of odd order | classification}
\label{square}

We are ready to formulate and prove our main result for odd square STGQs.
For the rest of this section, the GQ $\mS^x$ will be a thick STGQ of odd order $s$. By \S \ref{rootel}, $\mS^x$ is central.

%\subsection{Ghinelli's results}

\subsection{Structure of odd STGQs}

We deduce some structural results from \S \ref{rootel}. The following theorem is essentially contained in \cite{SDWKTHeis}, and uses the Payne-derived quadrangle $\mP(\mS^x,x)$.
We give a short different proof here using centrality.

\bl
\label{centraleqodds}
We have that $Z(H) \leq \mathbb{S}$, so that $\mathbb{S}$ is the center of $H$.
\el

{\em Proof}.\quad
Suppose $b \in Z(H) \setminus \mathbb{S}$. As $H = \bigcup_{A \in \mF}A^*$, we can write $b = \beta s$ with $\beta$ contained in some element $B$ of $\mF$, and $s \in \mathbb{S}$.
As $s \in Z(H)$, $\beta$ also must be in the center. This contradicts Corollary \ref{centrals}.
\eop \\

We sum up some structural results in the next theorem.

\begin{theorem}[Structure theorem]
\label{MTSTGQ}
Let $\mS = (\mS^{x},H)$ be an STGQ of order $q$, $q$ an odd prime power, and let $(\mF,\mF^*)$ be the $4$-gonal
family arising from (w.r.t. any point not collinear with $x$). The following properties are equivalent.
\begin{itemize}
\item[{\rm (i)}]
The members of $\mF^*$ are abelian.
\item[{\rm (ii)}]
$\mS^{x}$ satisfies $(M)_x$.
\item[{\rm (iii)}]
$H$ is a special group of exponent $p$, $q$ being a power of the prime $p$.
\item[{\rm (iv)}]
The group of symmetries about $x$ is contained in $Z(H)$.
\item[{\rm (v)}]
The group of symmetries about $x$ coincides with $Z(H)$.
\end{itemize}
\end{theorem}

{\em Proof}.\quad
Denote the group of symmetries about $x$ by $\mathbb{S}$ as usual.
As $x$ is a regular point, there is a projective plane $\pi_x$ of order $q$ associated to $x$. Let $y \sim x \ne y$ be any point; then $H_y$ induces an automorphism group
$H_y/\mathbb{S}$ on $\pi_x$ that clearly fixes all the points on $xy$ of the plane. So $H_y$ fixes any point of the line $xy$
of $\mS$. In fact $H/\mathbb{S}$ induces a translation group of the translation plane $\pi_x$, so $H/\mathbb{S}$ is elementary abelian.\\

Assume (i). Let $M \nI x$ be arbitrary, and let $m$ be so that $x \sim m \I M$. Then $H_m$ is abelian, and 
so $H_M$ is a normal subgroup of $H_m$. Since $\mathbb{S} \leq H_m$, it follows that $H_M$ fixes every line incident with $m$.
So (ii) is satisfied.\\

Assume (ii). It is well-known that $H$ is generated by all groups $H_N$, with $N \nI x$. 
Since (ii) holds, any commutator of the form $[h,s]$, with $h \in H_N$ and $s \in \mathbb{S}$, is trivial, so that 
$\mathbb{S}$ commutes with any group of the form $H_N$. Whence $\mathbb{S}$ is a subgroup of $Z(H)$, and (iv) is implied.\\

If (iv) holds, then (v) follows by Lemma \ref{centraleqodds}.\\

Suppose (iv) or (v) holds. Let $y \sim x \ne y$, and $M \I y$, $M \ne xy$. Then $H_y = \mathbb{S}H_M$. Since $H_M \cap \mathbb{S}$ is trivial
and $H/\mathbb{S}$ is abelian, $H_M \cong H_M/H_M \cap \mathbb{S}$ also is. 
So $H_y$ is abelian, and (i) follows.\\

We have shown that (i), (ii), (iv) and (v) are equivalent. We still have to show that any of these properties is equivalent
with (iii).
First suppose that $H$ is special; then $\vert Z(H)\vert = q$, and from Lemma \ref{centraleqodds} (iv) follows. 
Now suppose the conditions (i), (ii) and (iv) hold. Consider an element $\alpha \in H$; then we can write 

\begin{equation} \alpha = abc, \end{equation}

where $a \in A, b \in B$, $A$ and $B$ being members of the associated Kantor family, and $c$ a symmetry about $x$. 
By \S \ref{Exp}, we have that 

\begin{equation} (abc)^p = a^pb^p[b,a]^{p(p -1)/2}c^p,        \end{equation}

where $q$ is a power of the prime $p$. By (i) we know that $A$ and $B$ are elementary abelian. Also, by \cite[Theorem 2]{Fro}, see \S \ref{Exp},
$Z(H)$ is elementary abelian, so $\alpha^p = (abc)^p = \he$, and $H$ has exponent $p$. 
As $[H,H]H^p = \Phi(H)$,  it follows that $\Phi(H) = Z(H)$.
Now consider any two elements $\alpha, \beta \in H$; then $[\alpha,\beta]$ is a symmetry about $x$ (look at the action
on $\pi_x$), so $[H,H] \leq Z(H)$.  Let $A,B \in \mF$ be distinct, and choose $\alpha \in A^{\times}, \beta,\beta' \in B^{\times}$;
then $[\alpha,\beta] = [\alpha,\beta']$ if and only if $[\beta'\beta^{-1},\alpha] = \he$ if and only if $\beta'\beta^{-1}$
is a symmetry about $Z(H)$, cf. Corollary \ref{centrals}. This implies that $\beta = \beta'$, and one deduces that $\vert [H,H] \vert \geq q$, so
$[H,H] = Z(H)$.
\eop

\bl[Projection Lemma]
\label{Tproj}
No proper subgroup of $H$ projects onto $T$.
\el
{\em Proof}.\quad
By the previous lemma, we know that $\Phi(H) = \mathbb{S}$. Now, by way of contradiction, let $R \leq H \ne R$ project on $T$, and suppose $M \leq H$ is a
proper maximal subgroup containing $R$. Then $M$ contains $\Phi(H)$, so also $\mathbb{S}$; this leads to the fact that $M = H$, a contradiction. 
\eop

%\bigskip
\subsection{Property (AR1)}
Take $l \in A \in \mF$; then by Corollary \ref{centrals}

\begin{equation} [l,H] = [l,A^*B] = [l,C_H(l)B] = [l,B] = [H,H] = \mathbb{S}.\end{equation}
 
This implies  the following property which generalizes antiregularity:
 
 \begin{quote}
 (AR1)\ \ No triad $\{U, V, W\}$ of lines, with $U \I x$ and center $X \I x$, has more than one center which is not incident with $x$.
 \end{quote} 
 
 In the next section we will further exploit (AR1) to show that EGQs with (AR1) are classical.

\bigskip
\subsection{Main result}

A group $G$ is of {\em semifield type} if there is some prime power $q$ for which $\vert H\vert = q^3$, and such that there are 
normal elementary abelian subgroups $M$ and $N$ of order $q^2$ for which

\begin{itemize}
\item[(i)]
$H = MN$;
\item[(ii)]
$[m,n] = \he$ implies $m \in M\cap N$ or $n \in M \cap N$ for all $m \in M$ and $n \in N$.
\end{itemize} 

In \cite{Hira}, the author shows that from a semifield group $G$ and a pair $(M,N)$ such as above, a semifield $\mathbb{S}(M,N)$ can be constucted of order $q$, and conversely,
any semifield $\mathbb{S}$ gives rise to a semifield group in a natural way (if the plane is not Desarguesian, the group is simply  the collineation group generated by all elations). In fact, the corresponding semfield group acts naturally on the semifield plane $\Pi(\mathbb{S})$ which 
is coordinatized by $\mathbb{S}$.    \\

Let $A^*, B^*$ be distinct elements of $\mF^*$; then $A^*$ and $B^*$ are elementary abelian and normal (since they both contain $[H,H]$), and (ii) is obviously satisfied
by Corollary \ref{centrals}. So $H$ is a group semifield type with respect to any two distinct members of $\mF^*$.
We now address the following result of Hiramine:

\bt[Y. Hiramine \cite{Hira}]
\label{Hirathm}
Let $G$ be a group of semifield type of order $q^3$, with $q = p^h$, $p$ a prime and $h \in \mathbb{N}^\times$. Let $W_G = W$ be the set of all abelian subgroups of order $q^2$.  
\begin{itemize}
\item
If $p = 2$, then $\vert W \vert = 2$.
\item
If $p > 2$, then $\vert W\vert = p^r + 1$ for some natural number $r \leq h$, and $H \curvearrowright W$ contains the natural ($2$-transitive) action of $\PSL_2(p^r)$ of degree
$p^r + 1$.
\end{itemize}

Also, $\vert W \vert = q + 1$ if and only if the corresponding plane is Desarguesian; in that case, $G$ is isomorphic to the following group:

\begin{equation} \{\left(
 \begin{array}{ccc}
 1 & \alpha & c\\
 0 & \mathbb{I} & \beta^T\\
 0 & 0 & 1\\
 \end{array}
 \right)\vert  \alpha, \beta \in \mathbb{F}_{q^2}, c \in \mathbb{F}_q \} \cong \hH_1(q),                                         \end{equation}

with standard matrix multiplication, and where $\mathbb{I}$ is the identity $2 \times 2$-matrix.\\
\et

As $\vert \mF^*\vert = q + 1$, we have obtained the following theorem.

\bt
$H$ is isomorphic to the classical semifield group $\hH_1(q)$. \eop \\
\et

Now consider any line $[A] \I x$, $A \in \mF$, and consider $A^H$. By (AR1) (see also Corollary \ref{centrals}), any two of the $q$ members of $A^H$
meet only in the identity; it follows that $A^H \cup \mathbb{S}$ is a partition of $A^*$ into $q + 1$ subgroups of order $q$. Whence there arises a 
translation plane $\Pi_A$ of order $q$ for which $A^*$ is the translation group. Moreover, $H$ acts faithfully as a collineation group on $\Pi_A$, and since
for any $B \in \mF \setminus \{A\}$, $B$ fixes the point at infinity of $\Pi_A$ corresponding  to $\mathbb{S}$ and acts sharply transitively on the other points at infinity, 
$\Pi_A$ is a semifield plane. By Theorem \ref{Hirathm}, $\Pi_A$ is Desarguesian.

\br
{\rm Note that the planes $\Pi_A$ are nothing else than the planes we construct geometrically in the next section by using (AR1).}
\er

Consider the projective plane $\Pi_x$; it is the projective completion of the translation plane which arises from the partition in the group $H/\mathbb{S}$ defined
by $\{ A^*/\mathbb{S} \vert A^* \in \mF^* \}$. As $H \cong \hH_1(q)$ and as the number of maximal elementary abelian subgroups of $\hH_1(q)$ is $q + 1$, 
$\Pi_x$ must be Desarguesian. (Note that $[.,.]: H/\mathbb{S} \times H/\mathbb{S} \mapsto \mathbb{S}$ defines an alternating form over $\mathbb{F}_q$.)\\

At this point, there are different ways to proceed. One is to construct explicitly an automorphism of $H$ which sends the Kantor family of $\mS^x$ to the 
classical Kantor family of $\hW(t)$. As this is done in the recent paper \cite{Ghi}, we will not do this, and instead refer to {\em loc. cit.} | see \S \ref{Ghinsec}.
Another way is to ``extend (AR1)'' to the property that all lines incident with $x$ are antiregular (again, we won't do it explicitly due to Ghinelli's paper).  The result (Theorem \ref{mainsquareodd}) then follows from the fact that each plane $\Pi_A$ is Desarguesian and  \cite[5.2.7]{PTsec}.

%Now consider any triad $\{[A],U,V\}$, with $\proj_xU \ne \proj_xV$ and $A \in \mF$. If we show that such a triad has precisely $0$ or $2$ centers, (AR1) implies %that 
%each line incident with $x$ is antiregular (and so all lines of $\mS^x$ are). So suppose $W \in \{[A],U,V\}^{\perp}$, and suppose there are no other centers.
%Let $W \cap U = w$ and $W \cap V = w'$, and let $\proj_Ux = u$ and $\proj_Vx = v$; also, put $a = \proj_Wx$.  Each of the $q$ affine points of $U$, respectively %$V$,
%defines a line in the dual plane $\Pi_x^D$ (by considering the perp with $x$) through $u$, respectively $v$. Also, the lines $L(w)$ and $L(w')$ defined by 
%$w$ and $w'$ intersect in $a$. Consider the group $H_U$ (as a collineation group of $\Pi_x^D$), and the orbit $L(w)^{H_U}$ (these are just the lines of $\Pi_x^D$ %on $u$). For each such line $L(u')$, $u' \I U$, $u \ne u'$, define
%$p_{u'} = L(u') \cap L(\proj_Vu')$. Then using the action of $H_U$ and the fact that $\Pi_x$ is Desarguesian one observes that any line on $x$ except $\proj_xU$ %and 
%$\proj_xV$ contains $0$ or $2$ points  of $\{ p_{u'} \vert u' \I U, u' \not\sim x\}$ (compare the situation when $\mS^x = \hW(q)$). 
%So $\vert \{[A],U,V\}^{\perp}\vert = 2$, and indeed every line of $\mS^x$ is antiregular.

\bigskip
\bt
\label{mainsquareodd}
$\mS^x \cong \hW(q)$ and
$H \curvearrowright \mS^x$ is permutation equivalent with $\hH_1(q) \curvearrowright \hW(q)^z$ for any point $z \in \hW(q)$, the latter action being that as an elation group.\eop
\et

\bigskip
\subsection{Implications}

The next set of theorems follows immediately.

\bt
If $(\mS^x,H)$ is an EGQ of order $q$ with $H \cong \hH_1(q)$, then $\mS^x \cong \hW(q)$.
\eop \\
\et

\bt
If an EGQ of order $q$ admits an elation group of semifield type, we have the conclusion of Theorem \ref{mainsquareodd}.
\eop \\
\et

\bt
An EGQ with at least one antiregular line is isomorphic to $\hW(q)$.
\eop \\
\et

\bt
If $(\mS^x,H)$ is an EGQ of odd order $q$ and $\vert Z(H)\vert \geq q$, then we have the conclusion of Theorem \ref{mainsquareodd}.
\eop \\
\et

\bt
If $(\mS^x,H)$ is an EGQ of order $q$ with only abelian components (in the corresponding Kantor family), then we have the conclusion of Theorem \ref{mainsquareodd}.
\eop \\
\et

\medskip
\subsection{Ghinelli's result}
\label{Ghinsec}

In \cite{Ghi}, Ghinelli studies so-called ``AS-configurations'' in finite groups. Such a configuration is a set of subgroups $A_0,\ldots,A_{n + 1}$ in
a group $A$ of order $n^3$ ($n \in \mathbb{N}$, $n \geq 2$), such that $A_iA_j \cap A_k = \{ \id \}$ for two by two distinct $i, j, k$, while 
each $A_i$ has size $n$. Moreover, it is asked that $A_0 \unlhd A$. It is straightforward to observe that putting $\mF := \{ A_1,\ldots,A_{n + 1}\}$
and $\mF^* := \{ A_0A_1,\ldots,A_0A_{n + 1} \}$, one obtains a Kantor family of type $(n,n)$ in $A$, that is, an EGQ $\Gamma$ of order $n$ with elation 
group $A$. But there is more: the group $A_0$ obviously is a group of symmetries with center the elation point ($\mu$) | in other words, 
$\Gamma^{\mu}$ is an STGQ. Ghinelli then proves the following result (which we formulate in terms of STGQs, rather than in terms of AS-configurations).

\bt[Ghinelli \cite{Ghi}]
\label{Ghib}
A square STGQ $(\Gamma^{\mu},A)$ with odd order $n$ and of symplectic type is isomorphic to $\hW(n)$.
\et

{\em Symplectic type} means that each conjugacy class different from $\{ \id\}$ has a representative in $\cup_{j = 0}^{n + 1}A_j^{\times}$. 
This is nothing else than the group theoretical formulation of (AR1), which always holds for square STGQs of odd order.

\newpage
\section{Combinatorial lemma}
\label{comblem}

In this section, we observe that \cite[1.3.2]{PT} also holds for GQs that satisfy (AR1). We then derive some interesting consequences.\\

Let $Y$ and $X$ be concurrent lines in a GQ $\Gamma$ of order $s$, and suppose (AR1) is satisfied ``with respect to $X$ and $Y$''.
Put $X \cap Y = z$, and define a rank $2$ incidence structure $\Pi(X,Y) = (\mP,\mB,\I)$ as follows.

\begin{itemize}
\item
\textsc{Points}
are the elements of $X^{\perp}$ not through $z$.
\item
\textsc{Lines}
are of two types:
\begin{itemize}
\item
sets $\{Z,X\}^{\perp} \setminus \{Y\}$, with $Z \in Y^{\perp}$, $Z \nI z$;
\item
the points of $X \setminus \{z\}$.
\end{itemize}
\item
\textsc{Incidence}
is the expected one.\\
\end{itemize}

\bt
$\Pi(X,Y)$ is an affine plane of order $s$. 
\et
{\em Proof}.\quad
Let $y$ be any point incident with $Y$, and different from $z$. Clearly, if $V, W$ are such that $Y \ne V \I y \I W \ne Y$, then 

\begin{equation} \{V,X\}^{\perp} \cap \{W,X\}^{\perp} = \{Y\}. \ \ (\#) \end{equation} 

So
all lines in $\Pi(X,Y)$ determined by $y$ are parallel. 
Suppose $V, W \in Y^{\perp}$, $V \nI z \nI W$, and $V \not\sim W$. Then $\{V,X\}^{\perp}$ and $\{W,X\}^{\perp}$ meet in at most one line of $X^{\perp}$ not incident with $x$. Using 
$(\#)$ it follows that this is precisely one. The theorem follows.  
\eop \\

The set $X \setminus \{z\}$ is a parallel class of lines.
The other parallel classes correspond to points of $Y \setminus \{z\}$. Clearly, an element of $\mathrm{Aut}(\Gamma)_{Y,X}$ induces an element of $\mathrm{Aut}(\Pi(X,Y))$. \\

Below, if $Y \I z$ and the element of $\mF^*$ corresponding to $Y$ has (*), we also say that ``(*)$_Y$ is satisfied''. And if $(M)_x$ is true locally at $Y$ for $G$, we say that ``$(M)_{z,Y}$ is satisfied''.

\bt
Let $Y$ and $X$ be distinct concurrent lines in the EGQ $(\Gamma^z,G)$ of order $s$, $X \I z \I Y$, and suppose $(AR1)$ is satisfied with respect to $X$ and $Y$.
If $(\Gamma^z,G)$ satisfies (*)$_Y$, it also satisfies ${(M)_{z,Y}}$. 
\et
{\em Proof}.\quad
The subgroup $G_U$, $U$ any line in $Y^{\perp}$ not incident with $z$, induces a translation group (of order $s^2$) of $\Pi(X,Y)$. (Note that (*)$_Y$ is used.)
Whence $G_U$ fixes all lines on $\proj_Uz$, i.e., ${(M)_{z,Y}}$ is satisfied.
\eop \\

In particular, if $\Gamma^z$ satisfies (AR1) for all line pairs on $z$, and (*) holds, it also satisfies ${(M)_z}$.

\bc
Under the assumption of the previous theorem, we have that $G_y$ is elementary abelian for any $y \in Y\setminus \{z\}$. 
\eop \\
\ec

Now suppose again that $(\Gamma^z,G)$ is an EGQ of order $s$, with (AR1) for all line pairs on $z$.
Then 

\begin{equation}  G_U^\times \cap G_{U'}^\times = \emptyset       \end{equation} 

for lines $U, U'$ with $U \cap U' \not\sim z$. (Here non-concurrent lines are allowed!) So the only elements of $G$ which do not fix some affine line, are contained in 
a set of size $s - 1$. It immediately follows that $z$ is a center of symmetry, so that $\Gamma^z$ is an STGQ. (Whence (*) holds.)
The next corollary now follows from our analysis in the square case.

\bt
\label{AR1thm}
Let $\Gamma^z$ be an EGQ of order $s$.
If one of the following properties hold, then $\Gamma^z \cong \hW(s)$, $s$ is odd, and the elation group is isomorphic to $\hH_1(q)$.
\begin{itemize}
\item
Each pair of lines incident with $z$ satisfies $(AR1)$.
\item
Some line of $\Gamma^z$ is antiregular.\eop \\
\end{itemize}
\et

\br[Even case]
{\rm
Let $\Gamma^z$ be an EGQ, and suppose $X \I z \I Y$ are distinct regular  lines. 
Then it is easy to prove (for instance by applying the methods of \cite[Chapter 2]{SFGQ}) that $X$ and $Y$ are axes of symmetry. This is already sufficient
for $\Gamma^z$ to be a TGQ. When $s$ is even, $\Gamma^z$ is also an STGQ in this case. 
And of course, if some line not on $z$ is regular, $\Gamma^z \cong \mQ(4,s)$.
}
\er

Note the following property, which can be obtained ``directly'' from (AR1).

\bt
Let $(\mS^x,H)$ be an EGQ of order $s$, with $(AR1)$. Then $H$ has exponent $p$, $H$ being a $p$-group.
\et
{\em Proof}.\quad
By (AR1), $(M)_x$ is satisfied, and $H_U^\times \cap H_V^\times = \emptyset$ for $U \ne V$ nonconcurrent affine lines. 
So $\{\mathbb{S}\} \cup_{U \relbar\hspace*{-0.2cm}\I\hspace*{0.1cm} x} \{H_U\}$ defines an equal sized partition of $H$, implying that $H$ has exponent $p$. \eop \\

\quad\textsc{Question}.\quad
{\em Does (AR1) force the GQ to be square?}\\

%\bt
%If the EGQ $(\mS^x,G)$ of order $(s,t)$ has $(AR1)$, then $s = t$ is odd.
%\et
%{\em Proof}.\quad
%Clearly $t \geq s$. Counting the number $\theta$ of $G$-elements stabilizing a line (call the set of such elements $S$), we obtain 

%\begin{equation}  \theta \geq 1 + (t + 1)s(s - 1) = 1 + s^2t + s^2 - st - s.     \end{equation}

%Since $\vert G \setminus S\vert \leq s - 1$, and $G \setminus S$ contains $G_{y,z}^{\times}$ for all $y \not\sim z$ in $x^{\perp}$, it follows that $s = t$.
%\eop \\

\bc
If the EGQ $(\mS^x,G)$ of order $(s,t)$ has $(AR1)$ and $s = t$, then $t$ is odd and $G \curvearrowright \mS^x$ is permutation equivalent with $\hH_1(t) \curvearrowright \hW(t)^z$, with $z$ any point of $\hW(t)$.
\eop
\ec

\newpage
\vspace*{6cm}
\begin{center}
\item
{\bf PART IV}
\item
\item
{\bf STGQS WITH IDEAL SUBGQS}
\end{center}
\addcontentsline{toc}{chapter}{Part IV - STGQs with ideal subGQs}

\newpage

\section{Tame groups and twisting}

In this section we introduce the process of ``twisting''; it was first described in \cite{Basic}. To start, we need the definition of tame group, which, for now, we 
describe as follows.

\subsection{Tame groups}

Call an elation group of an STGQ $\mS^x$ {\em tame} if it satisfies $(M)_x$.
We conjecture that  ``most'' STGQs admit a tame elation group. Still, it could happen that 
other elation groups arise, as indicated in the examples section, and in this and the next section we will classify the STGQs
with this property in a very precise way. For this purpose, we introduce {\em twisting} \cite{Basic} | a way to construct other elation
groups  (``twists'') starting from data 
\begin{center}
(certain EGQ, involution fixing certain subGQ elementwise).
\end{center}

\subsection{Twisting}
{\em For now, $\mS^{x} = \mS$ is an EGQ of order $(q^2,q)$, $q$ even,  with elation group $H$.
Also, $\mS'$ is a subGQ of order $(s',q)$, $s' > 1$, which is
fixed elementwise by a nontrivial collineation $\theta$ of $\mS$.
By \cite{burnside}, we then have that $\mS' \cong \hW(q)$, and that $\theta$ is an involution.}\\

%\medskip
Suppose $W$ is the group of all whorls about $x$, and let $S_2$ be a Sylow $2$-subgroup of $W$ which contains $H$. Then $S_2$ clearly has size
$2q^5$. Put $H' = \theta H$, so that $S_2 = H \cup H'$.
As $\mS' \cong \hW(q)$, and as each point of $\hW(q)$ is regular, one observes that for each point $z \not\sim x$, the pair $\{x,z\}$ is regular, so that
$x$ is a regular point of $\mS$.
This implies that two distinct subGQs of order $q$ containing $p$ can only intersect in a very restricted
manner:
either they share the lines through $x$ and the points (of the subGQs) incident with these lines, or
they intersect in the points and lines of a dual grid of order $(1,q)$.
Let $\theta'$ and $\theta''$ be two distinct nontrivial involutions in $S_2$ that respectively fix
the subGQs $\mS_{\theta'}$ and $\mS_{\theta''}$ (of order $q$) pointwise. Suppose that they intersect in a
dual grid as above. Then there is a point $z \not\sim x$ for which $\{x,z\}^{\perp\perp} \subseteq
\mS_{\theta'}\cap\mS_{\theta''}$. Since both $\theta'$ and $\theta''$ fix $z$, we immediately have
a contradiction since $\theta' \ne \theta''$ and $\vert {(S_2)}_z\vert = 2$.

So subGQs of order $q$ that are fixed pointwise by a nontrivial involution in $S_2$
mutually do not share points not collinear with $x$. This implies that if $\mS_1$ and $\mS_2 \ne \mS_1$
are two such subGQs, there is some line $M \I p$ so that
\begin{equation} 
M^{\perp} \cap \mS_1 = M^{\perp} \cap \mS_2.
\end{equation}

Also, it follows easily that the number of such subGQs is $q^2$, and that the associated involutions
are mutually conjugate in $S_2$.
Note also that all whorls of $S_2$ which are not elations about $x$ are contained in $H'$.
The group $S_2$ is non-cyclic; if it were cyclic, then $H$ would be abelian, implying in its turn that there are more lines through a point
than points incident with a line (since $\mS$ is then a TGQ, cf. Chapter 8 of \cite{PT}).
As $S_2$ is non-cyclic, a result of P. Deligne \cite{Deligne} implies that $S_2$ has at least
three subgroups of size $q^5$ (one of which is $H$). Suppose $H'' \ne H$ is a subgroup of $S_2$ of order $q^5$. If $H''$ does not contain
any of the $q^2$ involutions of above, then $H''$ is an elation group (``first case''). If $H''$ contains at least one such involution (``second case''), it contains all of them since they are
mutually conjugate, and since $H''$ is a normal subgroup of $S_2$ (as a group of index $2$).
In that case, put $H_1 = H'' \cap H$, and $H_2 = H' \cap H''$. So $\vert H_1\vert = \vert H_2\vert = q^5/2$. Then it is straightforward to see that
\begin{equation}  
H_1 \cup \theta[H \setminus H_1] = H^-                    
\end{equation}
\noindent
is an elation group of size $q^5$.
As the first and second case are equivalent, we keep using the notation of the second case.
We put $H_4 = H \setminus H_1$ and $H_3 = \theta H_4$.\\

Suppose $L \I x$, and let $y \I L \I x \ne y$.
Define $H(x,L,y) := H_{[\{ x,L,y\}]}$.\\

{\em Second Hypothesis}.\quad
{\em For all $L \I x$ and $y \I L \I x \ne y$, we have that $\vert H(x,L,y)\vert = q^2$.}
{\em Also, $H^2 \leq Z(H)$, where $H^2 = \{h^2 \vert h \in H\}$ and $Z(H)$ is the center of $H$.}

\br{\rm
We leave it as an exercise to the reader to show that $x$ is a center of symmetry.}
\er

Since $\vert H(x,L,y)\vert = q^2$ for all $L$ and $y$ as above, and since these groups generate $H$,
it is straightforward to show that $Z(H)$ is the group of symmetries about $x$. In fact, one
observes that $Z(H) = Z(H^-)$.
Let $H(x,L,y)$ be a root-group; then $H(x,L,y)^2 \leq Z(H)$, so that $H(x,L,y)^2 = \{\id\}$.
So all such root-groups are elementary abelian. Now consider $\theta\phi \in H^-$, where $\phi \in
H_4$ is a non-trivial root-elation in $H(x,M,z)$ with $z \in \mS_{\theta}$ which does not
fix $\mS_{\theta}$ (it is an easy exercise that such a $\phi$ exists for suitable $z$). Then
$(\theta\phi)^2 = [\theta,\phi^{-1}] = [\theta,\phi]$
clearly cannot be the identity, while it fixes $z$ linewise. So $(\theta\phi)^2 \not\in Z(H^-)$, and hence
$H \not\cong H^-$. One can now conclude the next theorem.

\bt[Twisting, \cite{Basic}]
\label{FIX}
Let $\mS = (\mS^{x},H)$ be an EGQ of order $(q^2,q)$, where $q$ is even,
which contains a subGQ $\mS'$ of order $(s,q)$, $s > 1$, fixed pointwise by a nontrivial automorphism
$\theta$ of $\mS$. If $H^2 \leq Z(H)$, and if for all $L \I x$ and $y \I L \I x \ne y$, we have that $\vert H(x,L,y)\vert = q^2$,
then there is an automorphism group $H'$ of $\mS$ such that $H' \not\cong H$ and $(\mS^{x},H')$ is an EGQ.
\et

\bc[\cite{Basic}]
Let $\mS = (\mS^{x},H)$ be an EGQ of order $(q^2,q)$, where $q$ is even,
which contains a subGQ $\mS'$ of order $(s,q)$, $s > 1$, which is fixed pointwise by a nontrivial automorphism
$\theta$ of $\mS$.
Let $z \not\sim x$ and suppose $z \sim z_i \sim x$ for $i = 0,1,\ldots,q$. If all groups
$H(x,xz_i,z_i) \leq H$ are elementary abelian and have size $q^2$, then there is an automorphism
group $H'$ of $\mS$ such that $H' \not\cong H$ and $(\mS^{x},H')$ is an EGQ.
\ec

\newpage
\section{STGQs with more than one elation group are classical}
\label{morethan}

In this section it is our goal to give a complete answer to a question of Payne (asked at the 2004 Pingree Park conference ``Finite Geometries, Groups and Computation''). As the classical STGQ $\hH(3,t^2)$ is the only known finite quadrangle
which allows different (even nonisomorphic) elation groups, Payne asked whether this property characterized the quadrangle. We will show that the anser is ``yes''.

The (proof of the) next structure theorem is partly a prototype for the general classification of STGQs of order $(t^2,t)$ (even though it turns out that $\mathrm{char}(\mS^x) = 2$). 

\bt
\label{tametw}
Let $\mS^x$ be an STGQ of order $(s,t)$ such that $x$ is an elation point w.r.t. at least two distinct elation 
groups. Then 
\begin{itemize}
\item[{\rm (i)}]
$t = s^2$ is a power of $2$; 
\item[{\rm (ii)}]
$\mS^x$ contains $\hW(t)$-subGQs  through $x$ which are fixed pointwise by an involution;
\item[{\rm (iii)}]
any elation group w.r.t. $x$ is either tame or twisted;
\item[{\rm (iv)}]
there is a unique tame elation group $\widehat{E}$; furthermore, this group has the property that for each 
$y \in x^{\perp} \setminus \{x\}$, the i-root $(x,xy,y)$ is Moufang with root group in $\widehat{E}$, the latter being generated by all root-groups of the aforementioned type;
\item[{\rm (v)}]
any twisted elation group has class $3$, while the tame group has class $2$;
\item[{\rm (vi)}]
each elation group in (v) has exponent $4$.
\end{itemize}
\et

{\em Proof}.\quad
Suppose $G$ and $H$ are distinct elation groups for $\mS^x$, and let $W(x)$ be the group of all whorls about $x$, so that
$\vert W(x)\vert > \vert G\vert$. Let $p$ be the characteristic prime of $\mS^x$, and suppose $p$ divides $[W(x) : G]$.
Since $st$ is a $p$-power, a combination of Theorems \ref{8.1.1} and \ref{2.2.2} (recalling that $x$ is a regular point) yields
that any element of $W(x)_z$ of order $p^h$ ($h > 0$), $z$ being any opposite point to $x$, fixes a subGQ of order $(t,t)$ pointwise, while 
$s = t^2$. By Lemma 5.1 of \cite{Notes}, $p = 2$, and $[S(x) : G] = 2$, where $S(x)$ is a Sylow $2$-sugroup of $W(x)$ containing $G$.
Also, such a subGQ is isomorphic to $\hW(t)$.\\ 

Now suppose $p$ does not divide $[W(x) : G]$. Then all (``full'') elation groups about $x$ are Sylow $p$-subgroups in $W(x)$, so by Sylow's Theorem

\begin{equation} [W(x) : G] \geq p + 1.              \end{equation}

The Sylow subgroups are conjugate, so they all contain $\mathbb{S}$  (as usual the group of symmetries about $x$). Now consider the action of $W(x)$ on the 
set $T$ of spans of noncollinear points in $x^{\perp}$ | the kernel of this action is $\mathbb{S}$. If $\theta \in W(x)/\mathbb{S}$
fixes distinct elements of $T$, then either $\theta$ is trivial, or it fixes a subplane of order $t$ of the dual net $\mN^D_x$ (in which case $s = t^2$), cf. \cite{notenet}.\\
Suppose this latter case cannot occur {\em  for any} $\theta$; then $(W(x)/\mathbb{S},T)$ is a Frobenius group, and if $F$ is 
the Frobenius kernel, then $F'$, the unique subgroup of $W(x)$ containing $\mathbb{S}$ which projects $F$ onto
$W(x)/\mathbb{S}$, is a normal Sylow $p$-subgroup of $W(x)$, contradiction. Whence this situation does not occur.\\
Now assume some $\theta$ fixes a subplane of order $t$ elementwise, and choose $\theta$ in such a way that it is contained 
in $W(x)_z$, $z$ a point opposite $x$. (This choice can be done without loss of generality since $\mS^x$ is an EGQ.)
Then the fixed elements structure of $\theta$ is again a subGQ of order $t$, and $\theta$ is henceforth an involution in 
$\Aut(\mS^x)$. Whence by Lemma 5.1 of \cite{Notes}, $\mathrm{char}(\mS^x)$ is $2$, and $[W(x) : G] = 2$, contradicting our assumption.\\

This proves (i) and (ii).\\

The remaining cases are (obiously) less trivial. We proceed as follows.

First let $\{\theta_i \}$ be the set of all involutions in $W(x)$ that fix a $\hW(t)$-subGQ through $x$ pointwise, 
and define the injective map
\begin{equation} \Theta:\{\theta_i\} \hookrightarrow 2^{\mS^x}: \theta \mapsto \mS_{\theta},                 \end{equation}

which associates to each element of $\{\theta_i \}$  its fixed elements structure. Here, $2^{\mS^x}$ is the set of all subgeometries of $\mS^x$.
Any two of these involutions are clearly conjugate (in each Sylow $2$-subgroup of $W(x)$ such involutions 
are easily seen to be conjugate). Let $K$ be a Sylow $2$-subgroup of $W(x)$ containing some prescribed elation group 
 of $\mS^x$, and let $\{\theta_i \}_K$ denote $\{\theta_i \} \cap K$ | note that this set is not empty.
We are interested in how  distinct subGQs associated to elements of $\{\theta_i \}_K$ can intersect (which will
enable us to detect part of the structure of the STGQ).
Theoretically, ideal subGQs $\mS'$ and $\mS'' \ne \mS'$ could intersect in any of the following ways (cf. \cite[Lemma 4.2.5]{SFGQ} for a proof):

\begin{enumerate}
\item[$(1)$]
$\mathcal{S}' \cap \mathcal{S}''$ is a set of $t^2 + 1$ pairwise nonconcurrent lines
of $\mathcal{S}'$ {\em and} $\mathcal{S}''$;
\item[$(2)$]
$\mathcal{S}' \cap \mathcal{S}''$ consists of a line $L$ of $\mathcal{S}'$ (and $\mathcal{S}''$), together with all
points of $\mathcal{S}'$ (and $\mathcal{S}''$)
incident with this $L$, and all lines of $\mathcal{S}'$ (and $\mathcal{S}''$) incident with these points;
\item[$(3)$]
$\mathcal{S}' \cap \mathcal{S}''$  is a GQ of order $(1,t)$.
\end{enumerate}

Case (3) cannot occur by the twisting argument of the preceding section: since $\Theta$ is injective, 
$K_z$, $z$ being a hypothetical point of $\mS' \cap \mS''$ not in $x^{\perp}$, must contain distinct
involutions, contradicting its size. Case (1) can evidently not occur, so (2) remains.
It follows easily that 

\begin{equation} \vert \{\theta_i \}_K \vert = \vert \mbox{Im}(\Theta)_K \vert = t^2.            \end{equation}

(It is clear what we mean with $\mbox{Im}(\Theta)_K$.)\\
Now, following an idea of \cite{SolKnarr}, define an incidence structure $\Pi(K)$ as follows.

\begin{itemize}
\item
\textsc{Points} are the $t^2$ subGQs of $\mbox{Im}(\Theta)_K$ (``$K$-subGQs'').
\item
\textsc{Lines} are the point sets $\mS'' \cap M$, where $\mS'' \in \mbox{Im}(\Theta)_K$ and $M$ a line incident with $x$.
\item
\textsc{Incidence} is inverse containment.
\end{itemize}

Two distinct $K$-subGQs intersect in the $t + 1$ points on some line through $x$, together with all lines on these points. So in $\Pi(K)$, two distinct points are incident with exactly one line. Let $\mS''$ be a point, and $N$ be a line not containing $\mS''$ in $\Pi$.
Suppose $R \I x$ is the line which contains the point set $N$ in $\mS$. Then $\mS'' \cap R$ defines
the unique line of $\Pi$ parallel to $N$ and containing $\mS''$.
It follows that $\Pi$ is an affine plane of order $t$. 
Since $K$ induces a translation group on $\Pi(K)$, the following properties hold:

\begin{itemize}
\item[(a)]
if $k \in K$ fixes a $K$-subGQ, it fixes all $K$-subGQs;
\item[(b)]
by the previous property, $K_{\mS'}$ is a normal subgroup of $K$ for any $K$-subGQ $\mS'$, and 
$K/K_{\mS'}$ precisely is the translation group of $\Pi(K)$;
\item[(c)]
$[K,K] \leq K_{\mS'}$.
\end{itemize}

By (c) we have 
\begin{equation} [[K,K],[K,K]] \leq [K_{\mS'},K_{\mS'}] \leq K_{[\mS']},                       \end{equation}

the latter being a group of order $2$. (Remark that $\hW(t)$
is a TGQ since $\mathrm{char}(\mS^x) = 2$, so $K_{\mS'}/K_{[\mS']}$ is abelian.)
Whence if $H$ is an elation group w.r.t. $x$ in $K$, we have 
\begin{equation} [[H,H],[H,H]] \leq K_{[\mS']} \cap H = \{\he\}.                   \end{equation}

So the solvability length of $H$ is at most (and then precisely) $2$. We need to determine the nilpotency class.
Note first that since $[H,H] \leq H_{\mS'}$, $[H,H]$ is in the translation group of any of the $\hW(t)$'s earlier referred to, and in particular $[H,H]$ is abelian. 
Fix such a subGQ $\mS' \cong \hW(t)$, and let $T$ be its translation group w.r.t. $x$. Let $\gamma \in T^\times$; then  $\gamma$ fixes some line $L \I x$ pointwise in $\mS'$.
Since $H_{\mS'} \unlhd H$, $\gamma$ has this property in each element of ${\mS'}^H$, so $L$ is also fixed pointwise in $\mS^x$ by $\gamma$. So $T$ satisfies (*) as a subgroup of $\Aut(\mS^x)$.\\

Next we show that if 
 $U \nI x$, $U \in \mS'$, $H_U$ fixes $\proj_Ux$ linewise. For elements of $H_{\mS'} \cap H_U$ this is clear, so we consider an element $\zeta$ of $H_U$ not fixing $\mS'$. Then $H_U$ fixes the line $B$ of $\Pi(K)$ containing $\proj_Ux =:u$ (as a point set). As a point set, $B$ contains $t$ points of $xu$, so $H_U$
 fixes at least one point $u' \ne u$ in this set. In fact, $H_U$ fixes all points of $B$. To see this, let $\gamma \in H_U$ send $z$ to $z' \ne z$, both points in $B$.
 As $\gamma^2$ is a symmetry of $\mS'$ with axis $xu$, $\gamma$ acts as an involution on the affine lines which contain a point of $B$. Let $Z \I z$ be such an 
 affine line, and put $Z' := Z^{\gamma}$. Let $\upsilon$ be any element of $H_{\mS'}$ which sends $Z$ to $Z'$; it is also an involution, and fixes some line $F$
 incident with $x$ and different from $xu$ pointwise. It follows that $[\gamma,\upsilon] = \id$ (as $Z$ is fixed by $[\gamma,\upsilon]$, and $F$ pointwise). This means that $\gamma$ must fix the line $U^\upsilon \not\sim U$. It is easy to see that this contradicts the fact that $\{ U,U^{\upsilon}\}$ is a regular pair of lines in 
 each of the $K$-subGQs on $B$, together with the fact that any triad of lines in $\mS^x$ has precisley $t + 1$ centers. So indeed, $H_U$ fixes all points in $B$.
 
 Now take any $K$-subGQ on $B$, say $\widehat{\mS}$, and consider an affine line $C$ in $\widehat{\mS}$ which is not concurrent with $xu$. 
 Then  $[\zeta,H_C \cap H_{\mS'}] = \mathbb{S}$, for clearly $[\zeta,H_C \cap H_{\mS'}] \leq \mathbb{S}$, and $\zeta^a = \zeta^b$ with $a, b \in H_C \cap H_{\mS'}$ would imply that $\zeta^{ab} = \zeta$,
 which we have already seen not to be possible.  
 
 Note that for any $k \in H_C \cap H_{\mS'}$ we have
 \begin{equation}
 [\zeta,\zeta^k] = \zeta^{-2}[\zeta^{-1},k]\zeta^2[\zeta,k] =  [\zeta^{-1},k] [\zeta,k] = ([k,\zeta] [k,\zeta^{-1}])^{-1} = \id.
 \end{equation}
 (We use the fact that after the first equation symbol all elements are in $H_{\mS'}$, so commute with each other. Also, one notes that 
 $[k,\zeta]$ and $[k,\zeta^{-1}]$ are inverses of each other, as $\zeta$ and $\zeta^{-1}$ have the same action on $(xu)^{\perp} \cap \mS'$, and so both commutators yield the same symmetry in $\mathbb{S}$.)
 
 So indeed $\langle \zeta^k \vert k \in H_C \cap H_{\mS'} \rangle$ is an abelian group containing $\mathbb{S}$, which is what we wanted to prove.\\
 
 %Let $\mS'' = {\mS'}^{\zeta}$. 
 %Since $H_U$ has size $q^2 = 2^h$ for some positive integer $h$, $H_U$ must fix some line $U'$ on $\proj_Ux$ different from $U$ and not containing $x$.
 %Let $\eta \ne \zeta$ in $H_U$ also send $\mS'$ to $\mS''$; then $\zeta\eta^{-1}$ stabilizes $\mS'$, fixes $x$ linewise, and fixes both $U$ and $U'$. This implies %that $\zeta\eta^{-1}$ fixes 
 %$U \cap U'$ linewise (in $\mS'$, so also in $\mS^x$). (One easy way to see this is the following: $\zeta\eta^{-1}$ is an element of $\PGL_4(t)$ since it fixes $x$ %linewise, so
 %$U \cap U'$ must also be fixed linewise. But it is of course also just a symmetry in $\mS'$.)\\

Now consider $[[H,H],H]$; it is a subgroup of $T = H_{\mS'}$, so it  also satisfies (*) (as does $[T,H]$). 
We have a look at two types of commutators in $[T,H]$.

\begin{itemize}
\item
Type $[\alpha,\beta]$, with $\alpha \in T$ fixing $L \I x$ pointwise, and $\beta \in H$ fixing $y \sim x \ne y$ linewise, $y \nI L$.
Such a commutator is a symmetry about $x$.
\item
Type $[\alpha,\beta]$, with $\alpha \in T$ fixing $L \I x$ pointwise, and $\beta \in H$ fixing $y \sim x \ne y$ linewise, $y \I L$.
Then $[\alpha,\beta]$ is a root-elation in $\mS'$, so it is a symmetry about $L$ in $\mS'$. It easily follows that it also must be a symmetry about $L$ in $\mS^x$, which
contradicts $s = t^2$ and \cite[8.1.2]{PTsec}.
 \end{itemize}

Since $H$ is generated by the set $V$ of elements $\beta$ as above, we have that $[[H,H],H] \leq [T,H] \leq \mathbb{S}$.
Since again $H$ is generated by $V$, the nilpotency class of $H$ is thus at most $3$.\\

Class $1$ is not possible since $t < s$, an
obstruction for TGQs.\\

Now suppose the class of $H$ is precisely $2$, so that $\mathbb{S} \leq Z(H)$.
Then in $H$, all Moufang conditions are satisfied w.r.t. $x$ and $H$ is generated by the as such defined root groups (cf. \cite{SolKnarr}), so that 
$H$ is tame, and unique.  
(Note that, conversely, if $H$ would be tame, then clearly $\mathbb{S} \leq Z(H)$, and $H$ is of class $2$.)\\

Next suppose $H$ to have class $3$ | not all local Moufang conditions in $H$ w.r.t. $x$ are satisfied. 
We need to show that all i-roots and dual i-roots are Moufang, so that there indeed exists a tame elation group, which is then unique by the above.
So let $U \not\sim x$ be a line, and suppose $V = \proj_xU$. Let $S$ be a Sylow $2$-subgroup in $W(x)$. Construct the subGQ-plane $\Pi(S)$ w.r.t. $S$, and 
let $B$ be the line containing $V \cap U$ as a point set.  Then $S_U$ has order $2t^2$, and $S_U$ acts on $\Pi(S)$. Note the following facts:

\begin{itemize}
\item
$S_U$ fixes the point set corresponding to $B$ elementwise, and fixes each line incident with $V \cap U$;
\item
any element of $S_U$ which fixes a subGQ corresponding to a point incident with $B$ in $\Pi(S)$, either fixes it pointwise or induces a translation;
\item
if an element of $S_U$ fixes one, and then all, subGQs corresponding to $\Pi(S)$, then, since a translation in $\hW(t)$ always fixes some line pointwise,
it fixes $V$ pointwise if (and only if) it is an elation with center $x$. 
\end{itemize}

It follows that $\vert S_U/N\vert = 2$, with $N$ the kernel of  $S_U$ in its action on the points of $V$. So if $z \I V$, $z \not\in B \cup \{x\}$, then 
$S_{U,z}$ is a group of size $t^2$ which fixes each point incident with $V$. So $(V\cap U,V,x)$ is a Moufang i-root with root group $S_{U,z}$. This proves (iv)-(v).\\

(vi) follows easily from the fact that, if $\alpha \in H$, with $H$ some full elation group w.r.t. $x$, $H$ acting on the appropriate subGQ-plane (corresponding to a 
Sylow $2$-subgroup containing $H$), then $\alpha^2$ is the identity on the plane. So it fixes the subGQs in question, and hence it either is an involution fixing one 
of these subGQs pointwise (but then it is not in $H$), or it is a translation in one of these subGQs.
 \eop

\bc
$\mathbb{S} \leq Z(H)$ for any full elation group with center $x$.
\ec
{\em Proof}.\quad
Let $\widehat{E}$ be the tame elation group. Then $\widehat{E}$ is generated by the root-elations of $\mS^x$ with i-root containing $x$.
But each such elation commutes with any element of $\mathbb{S}$. So $\mathbb{S} \leq Z(\widehat{E})$. For non-tame elation groups, the argument is similar
(replacing root groups by groups of the form $H_U$, $U\not\sim x$, and remarking that $H_U$ fixes $\proj_Ux$ linewise).
\eop \\

\bc
\label{corfix}
Let $H$ be as in the theorem (a non-tame full elation group). Let $\gamma \in H \setminus \mathbb{S}$ fix points besides $x$. Then we have the following possibilities.
\begin{itemize}
\item[{\rm (i)}]
$\gamma$ fixes precisely $t + 1$ points, and $\gamma \not\in T$. Also, $\gamma$ fixes precisely one point $\ne x$ linewise.
\item[{\rm (ii)}]
$\gamma$ fixes $t^2 + 1$ points, and $\gamma \in T$. Also, $\gamma$ fixes $t^2$ affine lines.
\end{itemize}
\ec
{\em Proof}.\quad
By the previous theorem, we know that $H$ arises from twisting the unique tame elation group $L$ by some involution $\theta$ which fixes some 
ideal $\hW(t)$-subGQ elementwise. So $\gamma = \ell$ or $\gamma = \theta\ell$, with $\ell \in L$. This immediately gives us the 
information concerning the fixed points in (i) and (ii). If we are in case (i), then $\gamma \not\in T$ (as we have already seen), and by our first 
fixed point theorem, we also know that $\gamma$ must have fixed affine lines. The $t$ fixed points different from $x$ are all in some $\Pi(K)$-line (cf. the 
proof of the previous theorem), and by the same proof, all fixed lines must be incident with the same point. Clearly, this point is fixed linewise (by 
the fixed point theorem, or the fact that $\mathbb{S} \leq Z(H)$).

Now suppose that $\gamma$ fixes the $t^2 + 1$ points on some line incident with $x$.  By the fixed point theorem, $\gamma$ also fixes precisely $t^2$ affine lines.
We have seen that if $\gamma \not\in T$, all these lines must contain a point of one and the same $\Pi(K)$-line. But we have also seen that such $\gamma$ cannot exist.  
\eop \\

The next corollary is now immediate.
\bc
Let $H$ be as in the previous corollary. Then the i-root elations and dual i-root elations generate $T$.
\eop
\ec

We are ready to obtain the main result of this section.

\bt[Payne's question]
If $\mS^x$ is as in the previous theorem, it is isomorphic to $\hH(3,t^2)$.
\et

{\em Proof}.\quad
Let $H$ be as in the previous theorem. Let $\mS'$ be a $\hW(t)$-subGQ as above, and let $\mS''$ be another one in the same $H$-orbit.
Let $U \in \mS''$ be an affine line; it is exterior to $\mS'$, so it is concurrent with $t^2 + 1$ lines $T_U$ of $\mS'$ which partition its points (cf. the dual of Theorem \ref{2.2.1}). 
Define $A := H_U \cap H_{\mS'}$ | this subgroup of $H$ of size $t$ stabilizes $\mS'$, and also the  set $T_U$. As $A$ fixes $U$, Corollary \ref{corfix}
readily implies that $[U] := \proj_xU$ is fixed pointwise by $A$. Dualizing, and noting that $\hW(t)$ is self-dual as $t$ is even \cite[3.2.1]{PTsec}, $T_U$ becomes an ovoid $\mO$
of ${\mS'}^D \cong \hW(t)$ (which is also an ovoid in the ambient space of $\hW(t)$ | see \cite{ovoidal})
fixed by the group $A$; $A$ also fixes a flag $(X,[u])$ (corresponding to $(x,[U])$) ``elementwise'' (all the points on $X$ and 
all the lines on $u$). Let $\mO'$ be any $A$-orbit on $\mO$, different from $\{[u]\}$; each of its points is collinear with one and the same point of $X$, so 
the points of $\mO'$, together with $X$, generate a plane in the ambient $\PG(3,t)$ of $\hW(t)$. It easily follows that $\mO'$ is a translation oval, and so 
the $t$ $A$-orbits on $\mO$ different from $\{ [u]\}$ define a pencil of translation ovals with axis $X$, which lies in the tangent plane of $\mO$ at $[u]$.
By \cite{PentPraeg}, it follows that $\mO$ either is an elliptic quadric or a Suzuki-Tits ovoid. The main theorem of \cite{Brown} now implies that $\mS^x \cong \hH(3,t^2)$.
\eop \\

We will leave the next corollary (of the proof) to the reader.

\bc
Let $(\mS^x,H)$ be a finite EGQ of order $(s,t)$, and let $\mS'$ be a proper ideal subquadrangle of order $(t,t)$ containing $x$.
If $t$ is even and $\mS'$ is fixed pointwise by some involutory automorphism of $\mS^x$, $\mS^x \cong \hH(3,t^2)$.
\eop
\ec

(First show that $(\mS^x,H)$ is an STGQ, Then apply the previous theory.)

It might be in reach to conclude the same result as in the previous theorem, but with ``subGQ of order $(t,t)$ plus involution fixing it pointwise'' 
replaced by ``$\hW(t)$-subGQ''. I will come back to this issue in a future paper. \\

\newpage
\vspace*{6cm}
\begin{center}
\item
{\bf PART V}
\item
\item
{\bf GENERIC STGQs}
\end{center}
\addcontentsline{toc}{chapter}{Part V - Generic STGQs}

\newpage
\section{Generic STGQs}

When an STGQ $(\mS^x,H)$ satisfies (*), or equivalently when $T = H/\mathbb{S}$ is abelian, one easily observes that 
when $(\mF,\mF^*)$ is the associated Kantor-family, and $A, B$ are distinct members of $\mF$, that if $U \leq A$ and $V \leq B$, we have that
\begin{equation}
\langle U,V \rangle = H \ \ \longrightarrow \ \ U = A,\ \ V = B.
\end{equation}
(If at least one of $U, V$ is a proper subgroup, $\langle U\mathbb{S}/\mathbb{S}, V\mathbb{S}/\mathbb{S} \rangle$ is a proper 
subgroup of $T$.) We will show that a generalization of this natural property suffices to classify the parameters of the generic case (that is, the 
case where $T$ is not abelian), and to even classify to some extent. (So below, a {\em generic STGQ} is an STGQ for which $T$ is not abelian,
or equivalently, for which (*) is not satisfied. In particular, such STGQs cannot be square STGQs.)\\

More precisely, in a first instance our main objective is proving the following result. Before stating it, we need a definition. Let $(\mS^x,H)$ be as 
above, and let $\Phi := \Phi(H)$ be the Frattini subgroup of $H$. Define a point-line geometry $\Gamma(\Phi)$ as follows. 
Its lines are the $\Phi$-orbits on the lines incident with $x$ (where the trivial orbit $\{x\}$ is excluded); its points are the $\Phi$-orbits
in the affine point set of $\mS^x$. A point $u$ is incident with a line $V$ if at least one $\mS^x$-point of the orbit $V$ is collinear
with some point of the orbit $u$. If $u$ then is incident with $V$, it is easy to see that $V$ is surjectively projected on 
$u$, so that for each $L \I x$, a $\Gamma(\Phi)$-point is incident with precisely one line which is a $\Phi$-orbit on $L$. So each $\Gamma(\Phi)$-point
is incident with precisely $t + 1$ $\Gamma(\Phi)$-lines.

\bt
\label{Frattsub}
Let $(\mS^x,H)$ be a generic STGQ of order $(s,t)$, with Kantor family $(\mF,\mF^*)$. Suppose that the projection lemma (cf. the next
subsection) holds.
Suppose also that at least one of the following properties is satisfied.
\begin{itemize}
\item[(a)]
$\Gamma(\Phi)$ is a dual partial linear space.
\item[(b)]
For each $A \ne B \in \mF$, we have that if $K$ is a maximal subgroup of $H$ which does not contain $A$, then 
\begin{equation}
\langle A \cap K, B \rangle \ne H.
\end{equation} 
(This is the generalization of the aforementioned natural property.)
%\item[]
\end{itemize}
Then the Frattini subgroup $\Phi(H)$ is the elation group of a subSTGQ of order $t$, so that $s = t^2$.
\et

The way to construct ideal subquadrangles goes as follows.
Let $(\mK,\mK^*)$ be a Kantor family of type $(s,t)$ in the group $G$, and put $\hF := \mK \cup \mK^*$.
A nontrivial subgroup $X$ of $G$ is an {\em $\hF$-factor} of $G$ if
\begin{quote}
$(U\cap X)(V \cap X) = X$ for all $U,V \in \hF$ satisfying $UV = G$.
\end{quote}

Define $\mK_X = \{U\cap X\vert U \in \mK\}$ and $\mK^*_X = \{U^*\cap X \vert U^* \in \mK^*\}$.
We say that $X$ is ``of type $(\sigma,\tau)$'' if $\vert X\vert = \sigma^2\tau$, $\vert A\cap X\vert = \sigma$ and $\vert A^*\cap X\vert = \sigma\tau$ for all
$A \in \mK$ (and in \cite{Hach} it is shown that such integers $\sigma,\tau$ always exist).

\bt[D.~Hachenberger \cite{Hach}]
\label{Ha1}
Let $X$ be an $\hF$-factor of type $(\sigma,\tau)$ in $G$. Then necessarily one of the following cases occurs:
\begin{itemize}
\item[$(a)$]
$\sigma = 1$, $\vert X\vert = \tau \leq t$ and $X$ is a subgroup of $\cap_{A \in \mK}A^*$;
\item[$(b)$]
$\sigma > 1$, $\tau = t$ and $(\mK_X,\mK^*_X)$ is a Kantor family in $X$ of type $(\sigma,\tau)$.
\end{itemize}
\et

If we are in case (b) of Theorem \ref{Ha1}, we call $X$ a {\em thick} $\hF$-factor.
An $\hF$-factor $X$ in $G$ is {\em normal} if $X$ is a normal subgroup of $G$.

\bt[D. Hachenberger \cite{Hach}]
\label{Ha2}
Let $G$ be a group of order $s^2t$ admitting a Kantor family $(\mK,\mK^*)$ of type $(s,t)$, with $s,t > 1$, and having a normal $\hF$-factor $X$ of type $(\sigma,\tau)$ with $\tau = t$.
Then one of the following cases occurs:
\begin{itemize}
\item[$(a)$]
$G$ is a group of prime power order;
\item[$(b)$]
$\sigma > 1$, $\vert G\vert$ has exactly two prime divisors, and $X$ is a Sylow subgroup of $G$ for one of these primes.
\end{itemize}
\et

Once we have ideal subSTGQs, we can use our theory of square STGQs to search for more structure in the generic case.\\

The reader might want to keep the following standard lemma in mind:

\bl
\label{Frattder}
Let $G$ be a finite $p$-group. Then $\Phi(H) = [H,H]G^p$, where $G^p := \langle g^p \vert g \in G\rangle$.
\el

\br
{\rm In the course of the series, we will show that ``most of the time'', $H$ has exponent $p$, so that 
$\Phi(H) = [H,H]$. But when the exponent is $4$, it might happen that $\vert [H,H]\vert = t^3/2$ | in fact 
this happens for the exotic elation groups of $H(3,q^2)$ with $q$ even.}
\er

\subsection{Projection Lemma}

Let $(\mS^x,H)$ be an STGQ of order $(s,t)$. Let $T = H/\mathbb{S}$ be the translation group of the net $\Pi_x$.
In this paragraph we wish to classify the extensions

\begin{equation}  \he \mapsto M \cap \mathbb{S} \mapsto  M \mapsto T \mapsto \he,     \end{equation} 

with $M \leq H$. This will render us information about the Frattini subgroup of $H$.

The following lemma will be crucial (it is Theorem \ref{fix1}(iv) applied to STGQs with (*)).

\bl
\label{coll}
Suppose $\mS^x$ satisfies $(*)$.
If $s = t$ is odd, or $s \ne t$, any element $\alpha$ of $\cup_{A \in \hF}A^*$ maps some point of $\mP \setminus x^{\perp}$ to a collinear point.
\el
{\em Proof}.\quad
Suppose $\alpha \not\in \mathbb{S}$ does not have this property.
Applying Benson's lemma (while using $(*)$), we obtain

\begin{equation}  (t + 1)(s + 1)  + ts \equiv st + 1 \mod{s + t},    \end{equation}

so that $s + t$ divides $st$, contradiction. \eop \\

\bc
If $s = t$ is odd, or if $s \ne t$ and $T$ is abelian, the previous lemma applies,
\ec 
{\em Proof}.\quad
By assumption, $T$ is abelian, so $(*)$ holds for $\mS^x$ by Theorem \ref{ab1}. \eop \\

(As we will later see, $(*)$ will often be satisfied when $T$ is not abelian.)\\

Suppose $\mS^x$ has the property that any element $\alpha$ of $\cup_{A \in \hF} A^*$ maps some point of $\mP \setminus x^{\perp}$ to a collinear point. 
Suppose that $(*)$ holds.\\

Let $M$ be a subgroup as above, and note that without of loss of generality, we can suppose that $M$ is maximal in $H$. Let $(\hF,\hF^*)$ be the Kantor family
associated to $(\mS^x,H)$. Note that there must be some $A \in \hF$ not contained in $M$, since $\langle A \vert A \in \hF \rangle = H$.
Note that $[H : M] = p = [A^* : A^* \cap M]$, where $st$ is a power of the prime $p$. 
Let $L \I x$ be such that $L \ne [A]$, so that there are $p$ $(A \cap M)$-orbits on $L \setminus \{x\}$. On the other hand, $M$ acts transitively on $L \setminus \{x\}$
since $M$ projects on $T$. Moreover, since $[A^*: A^* \cap M] = p$ and $\mathbb{S} \not\leq M$, $A^* \cap M$ also acts transitively
on $L \setminus \{x\}$.
Let $\beta \in M \cap A^*$ be such that it does not fix all $(A \cap M)$-orbits on $L \setminus \{x\}$. Then $\beta$ fixes some line $U \not\in x$ (since 
it maps some affine point to a collinear point), and clearly $U \sim [A]$. Now $H_U$ is conjugate to $A$ by some element $\gamma$ of $H_{[L]}$, so 
$\beta$ is conjugate to some element of $A$ with the same action on $L$. But as $M$ is a normal subgroup of $H$, $\beta \in M$, so that 
$A \in M$, contradiction. 

\bt[Projection Lemma]
\label{ext}
The only extension of $T$ in $H$ is $H$ itself, if one of the following properties is satisfied.
\begin{itemize}
\item[{\rm (i)}]
$ s = t$ is odd.
\item[{\rm (ii)}]
$s \ne t$, and $\mS^x$ satisfies $(*)$. \eop
\end{itemize}
\et

%\medskip
\br[Even case]
{\rm
When $s = t$ is even, Theorem \ref{ext} cannot work in general. For, let $\mS^x$ be a TGQ of order $2^h$, so that $x$ is a regular point. 
Take any two distinct elements $A, B$ of the associated Kantor family. Then $AB = \langle A, B\rangle$, and $AB\mathbb{S}/\mathbb{S} = T$.
}
\er

 Let $M$ be maximal in $H$. By Theorem \ref{ext}, $M/(M \cap \mathbb{S})$ is an index $p$-group in $T$, so that $\mathbb{S} \leq M$.
 The following corollary is in some sense complementary to what we will encounter in the abelian factor section, later on. By $\Phi(H)$
we denote the Frattini subgroup of $H$.
 
 \bc[Structure of abelian $T$]
 Let $(\mS^x,H)$ satisfy one of $(i)$-$(ii)$ of the previous theorem. 
 \begin{itemize}
 \item[{\rm (i)}]
 If $A, B \in \hF$ are distinct, $\langle A, B \rangle = H$. Whence $\mathbb{S} \leq \Phi(H)$.
 \item[{\rm (ii)}]
 If $T$ is abelian, $\mathbb{S} = [H,H]$.
 \item[{\rm (iii)}]
 $\Phi(H) = \mathbb{S} = [H,H]$ if and only if $T$ is elementary abelian (so in particular, if $s = t$). 
 \item[{\rm (iv)}]
 If $T$ is abelian, every element of $\mathbb{S}$ is a commutator.
 \end{itemize}
 \ec
{\em Proof}.\quad
It is clear that $\langle A, B\rangle\mathbb{S}/\mathbb{S} = T$, so $\langle A, B\rangle = H$ by Theorem \ref{ext}. \\

Now let $T$ be abelian; then $[H,H] \leq \mathbb{S} = [H,H]$. Let $\gamma \in \mathbb{S}^\times$, and let $U \nI x$ be arbitrary.
Let $\xi$ be an element of $H \setminus \mathbb{S}$ mapping $U$ to $U^{\xi}$; then $\xi$ fixes some affine line $U' \ne U$ which meets 
$\proj_xU$. Let $\beta \in H$ map $U'$ to $U$. Then $[\beta,\xi]$ maps $U$ to $U^\xi$,  while being  a symmetry. So every element of $\mathbb{S}$ is a 
symmetry, and $[H,H] = \mathbb{S}$. Both groups coincide with $\Phi(H)$ if $T$ is elementary abelian. Vice versa, if $T$ is elementary abelian, 
\begin{equation} \bigcap_{Y \leq T\ \mbox{maximal}}Y = \{\he\},  \end{equation}
so that $\Phi(H) = \mathbb{S}$.  
\eop \\
 
 The following lemma follows in the same way as the Projection Lemma.
 
 \bt[Projection Lemma | general form]
\label{extg}
Suppose $(\mS^x,H)$ is an STGQ of order $(s,t)$, and let $(\hF,\hF^*)$ be the associated Kantor family.
Let one of the following properties be satisfied.
\begin{itemize}
\item[{\rm (i)}]
$ s = t$ is odd.
\item[{\rm (ii)}]
$s \ne t$, and $\mS^x$ satisfies $(*)$. 
\end{itemize}
Let $N$ be a normal subgroup of $H$, $L \I x$, and $O$ an $N$-orbit in $L \setminus \{x\}$. If $O$ is also an $(A^*\cap N)$-orbit for $A^* \in \hF^*$,
then $O$ is an $A$-orbit.     
\eop \\
\et

For a generic STGQ $(\mS^x,H)$, we say that {\em the projection lemma holds} if the following short exact sequence
\begin{equation}  \he \mapsto M \cap \mathbb{S} \mapsto  M \mapsto T \mapsto \he,     \end{equation} 
with $M \leq H$, only has solutions for $M = H$.

\bl
If the projection lemma holds for the generic STGQ $(\mS^x,H)$, then $\mathbb{S} \leq \Phi(H)$.
\el
{\em Proof.}\quad
Let $M$ be any maximal subgroup of $H$, then 
\begin{equation}
\frac{\vert H\vert}{\mathrm{char}(\mS^x)} = \vert M \vert = \vert M \cap \mathbb{S} \vert \times \vert M/(M \cap \mathbb{S}) \vert, 
\end{equation}
so that the projection lemma implies that $\mathbb{S} \leq M$. \eop \\

\subsection{Generic STGQs with (b)}

In this subsection we suppose that $(\mS^x,H)$ is a generic STGQ of order $(s,t)$ with Kantor family $(\mF,\mF^*)$, so that $s \ne t$. 
For the entire subsection, we assume to be in part (b) of the theorem. 
We may suppose w.l.o.g. that $\mathbb{S}$ is contained in $\Phi = \Phi(H)$.
Let $\mathbb{S} \leq G \leq H$, and $A,B \in \hF$ with $A \ne B$.  We say that $G$ satisfies (F) ``at $A$ w.r.t. $B^*$'' if 

\begin{equation}   G = (G \cap A)(G \cap B^*).    \end{equation}

Note that this property is symmetric in $A$ and $B$. If $G$ has (F) at $A$ w.r.t. any $B^* \in \hF^* \setminus \{A^*\}$, by definition $G$ satisfies (F) {\em at} $A$.\\

Let $A \in \hF$.
Define $\mR_A$ to be the set of maximal subgroups of $H$ that contain $A$; note that for $R \in \mR_A$, (F) is satisfied at $A$ w.r.t. any element of 
$\hF^* \setminus \{A^*\}$. Let $\mK_A := \cap_{C \in \mR_A}C$. Then the reader verifies that $\mK_A$ also satisfies (F) at $A$ w.r.t. any element of $\hF^* \setminus \{A^*\}$.
(If $A \leq L \leq H$, $L$ has (F) at $A$ w.r.t. any such element.)

\bl
\label{FAB}
Let $A, B$ be different elements in $\hF$.  If $K$ is a maximal subgroup in $H$ not containing $B$, then there is a maximal subgroup $K'$ containing 
$A$ for which $K' \cap B = K \cap B$. In particular, we have that $\mK_A \cap B \leq K \cap B \leq K$ for all maximal subgroups $K$. (So $\Phi \cap B = \mK_A \cap B$
if $A \ne B$ are elements in $\hF$, and whence for a fixed $B \in \hF$, $\mK_A \cap B$ is independent of the choice of $A \ne B$ in $\hF$.)
\el
{\em Proof}.\quad
As $\langle A,B\cap K\rangle \ne H$, there is a maximal subgroup $K'$ containing $\langle A,B\cap K\rangle$. Clearly, $K' \cap B = K \cap B$
(as otherwise $B \leq K'$, contradicting the projection lemma). \eop \\
%Consider the exact sequence

%\begin{equation} \he \mapsto \Phi \mapsto H \mapsto V \mapsto \he,         \end{equation}

%where we see  $V = H/\Phi$ as a vector space over $\mathbb{F}_p$. Suppose the hyperplane $K/\Phi$ intersects $B\Phi/\Phi$ in the space $\pi$.
%Then $\langle A\Phi/\Phi,\pi\rangle$ is a proper subspace of $V$, and its inverse image $\overline{\langle A\Phi/\Phi,\pi\rangle} =: K'$ in $H$ 
%contains $A$ and $B \cap K$, but not $B$.

%If $B \leq K$, we have that $\mK_A \cap B \leq K \cap B = B$; if not, then for some $K' \in \mR_A$ we have $K' \cap B = K \cap B$, so that $\mK_A \cap B \leq K' \cap B \leq K \cap B$.
%\eop \\

Lemma \ref{FAB} implies that for any $C, D \ne C \in \hF$, and any $A \in \hF$ (where $A = C$ or $A = D$ are allowed), 
\begin{equation}
A \cap \Phi = A \cap \mK_C \cap \mK_D.
\end{equation}
(Note that $\cap_{C \in \hF}C = \Phi$.) This is made slightly more explicit in the next lemma.

%\bl
%\label{FAB1}
%Let $\mathbb{S} \leq R \leq H$, and suppose $R$ satisfies $(F)$ at $A$ w.r.t. $B^*$. Let $K$ be a maximal subgroup of $H$, containing $B$. Then $K \cap R$
%satisfies $(F)$ at $A$ w.r.t. $B^*$. 
%\el

%{\em Proof}.\quad
%First note that $\mathbb{S} \leq K$, so $B^* \leq K$ while $[K :  K \cap A] = p$, where $st$ is a power of the prime $p$. If $R \leq K$, there is nothing 
%to prove, so suppose $[K : K \cap R] = p$. As $(K \cap R) \cap B^* = R \cap B^*$, the lemma follows from $R = (R \cap A)(R \cap B^*)$.
 %\eop \\

\bl
\label{FAB0}
Let $\mathbb{S} \leq L, M \leq H$, with $L$ and $M$ normal subgroups of $H$. Suppose that $L$ and $M$ have (F) at $A$ w.r.t. $B^*$ (and so also at $B$ w.r.t. $A^*$).
Then $L \cap M$ also has (F) at $A$ w.r.t. $B^*$. 
\el
{\em Proof}.\quad
We have that
\begin{equation}
\vert L \cap M \vert = \frac{\vert L\vert \times \vert M\vert}{\vert LM\vert} = \frac{(\vert A\cap L\vert \times \vert L \cap B^*\vert) \times (\vert A \cap M\vert \times \vert B^* \cap M \vert)}{\vert LM\vert}.
\end{equation}
As $\vert L \cap M \cap A \vert = \frac{\vert L \cap A \vert \times \vert A \cap M\vert}{\vert (L \cap A)(A \cap M)\vert}$ and 
$\vert L \cap M \cap B^* \vert = \frac{\vert L \cap B^* \vert \times \vert B^* \cap M\vert}{\vert (L \cap B^*)(B^* \cap M)\vert}$, it follows that
\begin{equation}
\label{eqpropf}
\vert L \cap M \vert =    (\vert L \cap M \cap A \vert \times \vert L \cap M \cap B^*\vert)\times\frac{\vert (L \cap A)(A \cap M)\vert\times\vert (L \cap B^*)(B^* \cap M)\vert}{\vert LM \vert}.
\end{equation}
Since $\vert L \cap M \vert \geq    \vert L \cap M \cap A \vert \times \vert L \cap M \cap B^*\vert$, we have that 
\begin{equation}
\frac{\vert (L \cap A)(A \cap M)\vert\times\vert (L \cap B^*)(B^* \cap M)\vert}{\vert LM \vert} \geq 1. 
\end{equation}
On the other hand, $(L \cap A)(A \cap M) \subseteq LM \cap A$ and $(L \cap B^*)(M \cap B^*) \subseteq LM \cap B^*$, while 
$(LM \cap A)(LM \cap B^*) \subseteq LM$. It follows that 
\begin{equation}
\frac{\vert (L \cap A)(A \cap M)\vert\times\vert (L \cap B^*)(B^* \cap M)\vert}{\vert LM \vert} =  1, 
\end{equation}
so that 
 $\vert L \cap M \vert =    \vert L \cap M \cap A \vert \times \vert L \cap M \cap B^*\vert$.
\eop

\bl
\label{FAB2}
Let $A, B$ ($\ne A$) $\in \hF$. Then $\mK_A \cap \mK_B$ satisfies $(F)$ at $A$ w.r.t. $B^*$. It follows that $\Phi = \mK_A \cap \mK_B$ for any choice of different $A, B \in \hF$.
\el

{\em Proof}.\quad
Consider $\mK_{A} \cap \mK_{B}$. Then applying Lemma \ref{FAB0} with $L = \mK_{A}$ and $M = \mK_B$, we conclude that $\mK_{A} \cap \mK_{B}$ has (F) at $A$ w.r.t. $B^*$. 
 
Noting that $\Phi \leq \mK_A \cap \mK_B = (\mK_A \cap \mK_B \cap A)(\mK_A \cap \mK_B \cap B^*)$, the fact that  $\mK_A \cap \mK_B \cap A = \Phi \cap A \leq \Phi$ and
$\mK_A \cap \mK_B \cap B^* = \Phi \cap B^* \leq \Phi$, indeed implies that $\Phi = \mK_A \cap \mK_B$ for any choice of different $A, B \in \hF$.
%Let $C \in \hF \setminus \{A_0,A_1\}$.
%Then $\mK_{A_1} \cap C \leq K$ for $K \in \mR_{A_2}$ by Lemma \ref{FAB}; whence $\mK_{A_1} \cap \mK_{A_2} \cap C = \mK_{A_1} \cap C$ and 
 %$\mK_{A_1} \cap \mK_{A_2} \cap A_1^* = \mK_{A_2} \cap A_1^*$. Now 
 %\begin{equation}
 %\vert \mK_{A_1} \cap C \vert \cdot \vert \mK_{A_2} \cap A_1^* \vert = \frac{\vert \mK_{A_2} \vert \cdot \vert C\vert}{\vert H\vert}\times\frac{\vert\mK_{A_2}\vert \cdot \vert A_1^*\vert}%{\vert H\vert} = \vert \mK_{A_1} \cap \mK_{A_2}\vert.
 %\end{equation} 
%So  $\mK_{A_1} \cap \mK_{A_2}$, has (F) at $A_1$ and $A_2$.
%The lemma easily follows by induction, and Lemma \ref{FAB0}.
\eop \\

\bc
Under the assumptions of this subsection, $\Phi = \Phi(H)$ is a normal $\hF$-factor in $H$. 
\eop
\ec

\medskip
\subsection{Generic STGQs with (a)}

Now we suppose that $(\mS^x,H)$ is a generic STGQ of order $(s,t)$ with Kantor family $(\mF,\mF^*)$,  
and  we assume to be in part (a) of the theorem. Let $A \ne B \in \mF$, and suppose that $abs = \phi$ with 
$a \in A$, $b \in B$, $s \in \mathbb{S}$, and $\varphi \in \Phi$. (Note that such $a, b, s$ always exist.) Then $ab = \phi s^{-1} =: \varphi \in \Phi$.
As $ab \in \Phi$, $ab$ fixes 
all points and lines of $\Gamma(\Phi)$. As $\Phi \unlhd H$, $a$ fixes all lines of $\Gamma(\Phi)$ of which the corresponding 
orbits are subsets of $[B]$ (``$B$-lines''), and $b$ fixes all lines of $\Gamma(\Phi)$ of which the corresponding 
orbits are subsets of $[A]$ (``$A$-lines''). As $\Gamma(\Phi)$ is a dual partial linear space, and as each $\Gamma(\Phi)$-point is 
incident with one $A$-line and one $B$-line, it follows that $a$ and $b$ fix all $\Gamma(\Phi)$-points (and so all $\Gamma(\Phi)$-lines).
One immediately deduces that $a \in \Phi$ and $b \in \Phi$.  So for all $C \ne D \in \mF$, we have shown that 
\begin{equation}
\Phi = (C \cap \Phi)(D^* \cap \Phi) = (D \cap \Phi)(C^* \cap \Phi).
\end{equation}

\bp
Under the assumptions of this subsection, $\Phi$ is a normal $\mF$-factor for $H$.
\eop \\
\ep

\subsection{Proof of Theorem \ref{Frattsub}}

 We have shown, under the assumptions of Theorem \ref{Frattsub}, that $\hF_{\Phi} := \{\Phi \cap A \vert A \in \hF\}$ defines a Kantor family of type $(s',t)$
in $\Phi$, with $s > s' > 1$, and $\hF^*_{\Phi} = \{\Phi \cap A^* \vert A \in \hF\}$. Also, since $x$ is regular, we have that $s'  = t = \sqrt{s}$ by Theorem \ref{2.2.2}.\\

\newpage
\vspace*{6cm}
\begin{center}
\item
{\bf APPENDIX}
\item
\item
{\bf REMARKS ON ABELIAN $T$, LOCAL (*) AND SUBGQ PLANES}
\end{center}
\addcontentsline{toc}{chapter}{Appendix - Remarks on abelian $T$, local (*) and subGQ planes}

\newpage
\section{Abelian factors, local (*) property and subGQ planes}

When $s = t$, we have seen that each member of the associated Kantor family $\hF$ is abelian. 
In this appendix, we first want to consider STGQs $(\mS^x,H)$  of general order $(s,t)$ with this property.
We will indicate that understanding a more general version of the subGQ plane which we encountered in \S \ref{morethan}, will 
be quite necessary in order to understand this situation, as well as the generic case.
Also, we will make some remarks on the class of STGQs which satisfy a (purely) local version of (*). 

These considerations serve mostly as an intro
to part II.

\subsection{Some results on PCPs}

Under the assumption that all elements of $\hF$ are abelian, the associated net $\Pi_x$, or rather associated ``PCP'' (see below), in $H/\mathbb{S}$ only consists of abelian factor groups $A\mathbb{S}/\mathbb{S}$, $A \in \hF$, so it makes sense to recall some basic facts about PCPs with
only abelian components.

Let $G$ be a group of order $s^2$, and let $\{G_1,\ldots,G_r\}$ $= \mG$ be a set of $ r \geq 3$ subgroups satisfying
\begin{itemize}
\item[(i)]
$\vert G_i\vert = s$ for $i = 1,\ldots,r$;
\item[(ii)]
$G_iG_j = G$ for distinct $i, j$.
\end{itemize}

Then $\mG$ is called a {\em partial congruence partition} (PCP) in $G$ with parameters $s$ and $r$, or a ``$(s,r)$-PCP'', with components the $H_i$s. The incidence structure
$(G,\{G_ig\})$ is a translation net with parameters $(s,r)$ and translation group $G$ (and vice versa). 

\bt[A. Sprague \cite{Sprague}]
\label{Spr}
If $G_1$ is normal in $G$, $G_2 \cong G_3 \cong \cdots \cong G_r$. If $G_1$ and $G_2$ are normal, all components are isomorphic, and $G \cong G_1 \times G_2$.
\et

As useful tool is the principle of {\em factorization}.

\begin{quote}
\textsc{Factorization}.\quad
Let $G$ and $\mG$ be as above. If $N$ is a normal subgroup of $G$ for which $N = (H \cap N)(K \cap N)$ for every distinct $H, K \in \mG$, then $N$ ``factorizes'' $\mG$.
We briefly say that $N$ has (F) in $\mG$. 
\end{quote}

\bo
If $N$ is as such, $G$ induces an $(n,r)$-PCP in $N$, where $\vert N\vert = n^2$; also, in $G/N$ an $(s/n,r)$-PCP arises.\\
\eo

Now suppose $(\mS^x,H)$ is an STGQ of order $(s,t)$, $t \ne s$, and suppose all elements of the associated Kantor family $\hF$ are abelian.
Then $\{A\mathbb{S}/\mathbb{S}\vert A \in \hF\}$  is an $(s,t + 1)$-PCP in  $H/\mathbb{S}$. (The PCP just corresponds to the net which arises from the 
regular point $x$.) 
By \cite{HachAI}, we want to consider two cases (which we will call ``AI-factors'' and ``AII-factors''): 

\begin{itemize}
\item
\textsc{AI-Factors.}\quad
$H/\mathbb{S}$ is abelian.
\item
\textsc{AII-Factors.}\quad
$H/\mathbb{S}$ has class $2$. In this case $Z(H/\mathbb{S})$ ($\ne \{\he\}$) factorizes $H/\mathbb{S}$.
\end{itemize}

Only those cases can arise, as will  be clear from the rest of this section. (The class of $H/\mathbb{S}$ determines the length of a tower of ideal subGQs | 
that is why AI-factors will be much harder to handle than AII-factors.)
The following theorem of Hachenberger uses the STGQ-PCP connection, and considers EGQs with {\em abelian Kantor families} (which are defined by the property that all members of $\mF \cup \mF^*$ are abelian). It essentially represents a special case of the one we are considering in the 
present section. We will show eventually that STGQs with abelian factors have abelian Kantor families. 

\bt[D. Hachenberger \cite{Hach}]
\label{HachAb}
Let $G$ be a nonabelian group admitting an abelian Kantor family $(\hF,\hF^*)$ of type $(s,t)$. Put $\Gamma = \hF \cup \hF^*$ and let $Z$ be the center of $G$.
Then necessarily $G$ is a group of prime power order, and one of the following occurs:
\begin{itemize}
\item[(a)]
$G/Z$ is elementary abelian;
\item[(b)]
$G$ has class $3$, and the members of the upper central series of $G$ are normal $\Gamma$-factors. The quadrangle is a central STGQ with $\mathbb{S} = Z$, and $t = s^2$;
\item[(c)]
$G$ is a $2$-group of class $2$, $G/Z$ is abelian but not elementary abelian, and the subgroup $Y = \langle g\vert (gZ)^2 = Z \rangle$ is a normal $\Gamma$-factor properly containing $Z$. We have that $t = s^2$, and the quadrangle is a central STGQ with $\mathbb{S} = Z$. 
\end{itemize}
\et

\medskip
\subsection{AI-Factors}

In this subsection $H/\mathbb{S}$ is abelian, and $p$ is the characteristic prime.
Since $H/\mathbb{S}$ is abelian, $[H,H] \leq \mathbb{S}$. (In particular, if $s = t^2$, $[H,H] = \mathbb{S}$.)\\

Suppose that $H/\mathbb{S}$ is {\em not} elementary abelian. Define 
\begin{equation} \Omega_p = \langle g\mathbb{S} \vert g \in H, (g\mathbb{S})^p = g\mathbb{S} \rangle \leq H/\mathbb{S}.\end{equation} 
Then by \cite{HachAI}, $\Omega_p$ factorizes in $H/\mathbb{S}$, and constitutes a $(t,t + 1)$-PCP, so a plane of order $t$ in $\Pi_x$.
Also, $p = 2$, and we are in Theorem \ref{HachAb}, third case.

\medskip
\subsection{AII-Factors}

Suppose $H/\mathbb{S}$ has class $2$.
For AII-factors, $Z(H/\mathbb{S})$ factorizes $H/\mathbb{S}$, so consitutes a $(t,t + 1)$-PCP in $H/\mathbb{S}$. Whence the net $\Pi_x$ contains 
a subplane of order $t$, with translation group $Z(H/\mathbb{S})$, and $s = t^2$. By Theorem \ref{netten}, $\mS^x$ contains
an ideal subGQ $\mS'$ of order $t$, and if $K \leq H$ is such that $K/\mathbb{S} = Z(H/\mathbb{S})$, then $K$ is an elation group of $\mS'$, forcing $\mS'$ to be an STGQ since $\mathbb{S} \leq K$. Note that $K$ is normal in $H$. \\

\medskip
\subsection{SubGQ plane}

Let $(\mS^x,H)$ be an STGQ. Let $\mS'$ be an ideal subGQ of order $t$ | so $\mS^x$ is of order $(t^2,t)$ | such that the following conditions
are satisfied:
\begin{itemize}
\item[(SUB-A)]
$({\mS'}^x,K)$ is an STGQ, for some $K \leq H$;
\item[(SUB-B)]
From $\mS'$ arises a subGQ plane (of order $t$), as described in \S \ref{morethan} | this amounts to asking that $K$ is normal in $H$, so that
${\mS'}^H$ consists of $t^2$ subGQs which two by two intersect in $t + 1$ points of some line on $x$, and the lines of $\mS$ incident with these points.
\end{itemize}

From the analysis of section \ref{morethan}, one can extract a great deal of properties from this situation, taken that $t$ is even.
Without assumptions on the parity of $t$, one can still prove that $\widehat{(M)_x}$ is satisfied, so that $\mS^x$ is central. 
More details will appear further in the series.\\

Both for AI-factors (taken that $H/\mathbb{S}$ is not elementary abelian) and AII-factors, such a subGQ plane arises. Also, we have shown
that all members in the class of generic STGQs contain subSTGQs of order $(t,t)$ which indeed satisfy (SUB-A) and (SUB-B) for the 
elation group $\Phi$. 

Again, this asks for the need to further analyze STGQs with subGQ planes.

\medskip
\subsection{STGQs with (*) at some element of $\mF^*$}

In Observation \ref{Structob}, we showed the following properties for the dual Suzuki-Tits quadrangles $\mS^x \cong \Gamma^D$ (we freely use the notation 
of that section):

\begin{itemize}
\item[{\rm (i)}]
$Z(G)$ is the symmetry group w.r.t. the elation point, and is elementary abelian.
\item[{\rm (ii)}]
$A^*(\infty)$ and $A(\infty)$ are elementary abelian and $A^*(\infty) \unlhd G$, so that $A^*(\infty)$ fixes $[A(\infty)]$ pointwise.
\item[{\rm (iii)}]
For $t \in \mathbb{F}_q$, $A^*(t)$ and $A(t)$ are non-abelian of exponent $4$; moreover, for $t \ne t'$, $A^*(t) \cong A^*(t')$ and $A(t) \cong A(t')$.
Also, no $A^*(t)$ is normal in $G$ and hence does not fix $[A(t)]$ pointwise.
\item[{\rm (iv)}]
$G$ is nonabelian of exponent $4$ and $G/Z(G)$ is not abelian.
\item[{\rm (v)}]
$G$ is the complete set of elations about $x$.
\item[{\rm (vi)}]
$\Gamma^D$ is not an MSTGQ.
\item[{\rm (vii)}]
$[G,G]$ is strictly contained in $A^*(\infty)$.\\
\end{itemize}

In part II, we will study STGQs which satisfy (*) at one member $A^* \leq G$ of $\mF^*$. (Note that by Theorem \ref{Spr}, if an STGQ satisfies (*) at at least two different members of $\mF^*$, all elements of $\mF$ are isomorphic.) 
We conjecture the following:\\

\quad\textsc{Conjecture} | [``Zero, one or all'']\quad
{\em For an STGQ $(\mS^x,K)$ with Kantor family $(\mF,\mF^*)$, either no, one or all elements of $\mF^*$ satisfy (*).}

%\newpage
%\section{STGQ parameter conjecture(s); Kantor conjecture}

\end{document}